\documentclass[twoside]{amsart}

\usepackage[latin1]{inputenc}

\newcommand{\SFr}{{\rm SFr}}
\newcommand{\SO}{{\rm SO}}
\newcommand{\R}{{\bf R}}

\newcommand{\Q}{{\bf Q}}

\newcommand{\Z}{{\bf Z}}
\newcommand{\C}{{\bf C}}
\newcommand{\CP}{{\bf CP}}

\newcommand{\Aa}{{\mathcal A}}
\newcommand{\Cc}{{\mathcal C}}
\newcommand{\Ee}{{\mathcal E}}

\newcommand{\Mm}{{\mathcal M}}

\newcommand{\Jj}{{\mathcal J}}
\newcommand{\Ss}{{\mathcal S}}
\newcommand{\Ll}{{\mathcal L}}
\newcommand{\Uu}{{\mathcal U}}

\newcommand{\Nn}{{\mathcal N}}
\newcommand{\id}{{\rm id}}
\newcommand{{\TS}}{{\Tilde {\rm Symp}}}
\newcommand{\p}{{\partial}}
\newcommand{\be}{{\beta}}
\newcommand{\Om}{{\Omega}}
\newcommand{\om}{{\omega}}
\newcommand{\eps}{{\varepsilon}}%
\newcommand{\de}{{\delta}}

\newcommand{\la}{{\lambda}}
\newcommand{\La}{{\Lambda}}
\newcommand{\si}{{\sigma}}

\newcommand{\Emb}{{\rm Emb}}

\newcommand{\Symp}{{\rm Symp}}
\newcommand{\Diff}{{\rm Diff}}

\newcommand{\Map}{{\rm Map}}

\newcommand{\Aut}{{\rm Aut}}

\newcommand{\PD}{{\rm PD}}

\newcommand{\U}{\mathrm{U}}
\newcommand{\PSU}{\mathrm{PSU}}
\newcommand{\Si}{{\Sigma}}

\newcommand{\MS}{{\medskip}}
\newcommand{\BS}{{\bigskip}}

\newcommand{\NI}{{\noindent}}
\newtheorem{theorem}{Theorem}[section]
\newtheorem*{theorem*}{Theorem}
\newtheorem{thm}[theorem]{Theorem}
\newtheorem{cor}[theorem]{Corollary}
\newtheorem{lemma}[theorem]{Lemma}
\newtheorem{definition}[theorem]{Definition}

\newtheorem{prop}[theorem]{Proposition}

\newtheorem{remark}[theorem]{Remark}

\newcommand{\omlc}{\tilde{\omega}_{\mu,c}}
\newcommand{\tH}{\tilde{H}}
\newcommand{\tT}{\tilde{T}}
\newcommand{\tJ}{\tilde{J}}
\newcommand{\jj}{{\mathcal J}}
\newcommand{\tjj}{\tilde{{\mathcal J}}}

\newcommand{\tM}{\tilde{M}}
\newcommand{\tMuc}{\tilde{M}_{\mu,c}}

\newcommand{\qq}{\mbox{\bf Q}}

\begin{document}

\title[Symplectic balls in rational $4$-manifolds]{The topology of the space of symplectic balls in rational $4$-manifolds}
\author{Fran\c{c}ois Lalonde}
\address{Universit\'e de Montr\'eal}
\email{lalonde@dms.umontreal.ca}
\author{Martin Pinsonnault}
\address{Universit\'e du Qu\'ebec \`a Montr\'eal}
\curraddr{University of Toronto}
\email{pinso@math.toronto.edu}
\date{May 31 2004}
\thanks{Lalonde's work is partially supported by NSERC OGP 0092913 and FCAR ER-1199}

\begin{abstract}
We study in this paper the rational homotopy type of the space of symplectic embeddings of the standard ball $B^4(c) \subset \R^4$ into $4$-dimensional rational symplectic manifolds, where $c= \pi r^2$ is the capacity of the standard ball of radius $r$. We compute the rational homotopy groups of that space when the $4$-manifold has the form $M_{\mu}= (S^2 \times S^2, \mu \om_0 \oplus \om_0)$, where $\om_0$ is the area form on the sphere with total area $1$ and $\mu$ belongs to the interval $[1,2]$. We show that, when $\mu$ is $1$, this space retracts to the space of symplectic frames for any value of $c$. However, for any given $1 < \mu \le 2$, the rational homotopy type of that space changes as $c$ crosses the critical parameter $\lambda = \mu - 1$, which is the difference of areas between the two $S^2$-factors. We prove, moreover, that the full homotopy type of that space changes only at that value, that is, that the restriction map between these spaces  is a homotopy equivalence as long as these values of $c$ remain either below or above that critical value.  The same methods apply to all other values of $\mu$ and other rational $4$-manifolds as well. The methods rely on two different tools:  the study of the action of symplectic groups on the stratified space of almost complex structures developed by Gromov, Abreu and McDuff and the analysis of the relations between the group corresponding to a manifold $M$, the group corresponding to its blow-up $\tM$, and the space of symplectic embedded balls in $M$.
\end{abstract}

\maketitle

\tableofcontents

\section{Introduction} \label{se:intro}

It is known that any symplectic form on the manifold $M=S^2 \times S^2$ is, up to a constant, diffeotopic to the split structure $M_{\mu}=(M,\om_{\mu})=(S^2\times S^2,\mu\om(1)\oplus\om(1))$, where $\mu$ is any real number belonging to $[1,\infty)$ and $\om(1)$ is the area form of total area equal to $1$ (see, for instance, Lalonde-McDuff~\cite{LMD}). Let $B^4(c)\subset\R^4$ be the closed standard ball of radius $r$ and capacity $c=\pi r^2$ equipped with the restriction of the symplectic structure $\om_{st}=dx_1\wedge dy_1+dx_2\wedge dy_2$ of $\R^4$. Let $\Emb_{\om}(c, \mu)$ be the space, endowed with the $C^{\infty}$-topology, of all symplectic embeddings of $B^4(c)$ in $M_{\mu}$. By the nonsqueezing theorem, this set is nonempty if and only if $c < 1$. McDuff~\cite{MD:BlowUps} showed that it is always connected.

In this paper, we study the homotopy type of that space for any given value $\mu\in [1,2]$ and all values of $0<c<1$. The general case $\mu\in [1,\infty)$ will be treated in~\cite{Pinso}, it involves more elaborated algebraic computations but the main geometric ideas are already present in the case $\mu\le 2$ of this paper. It will be convenient to denote by $[\mu]_-$ the lower integral part of $\mu$, that is, the largest integer amongst all those that are strictly smaller than $\mu$, and by $\la$ the difference $\mu-[\mu]_-$. Thus when $\mu$ belongs to $(1,2]$, $\la$ equals $\mu -1$, the difference of areas between the base and the fiber of $M_{\mu}$.

Consider the fibration
\begin{equation}\label{Fibration-1}
\Symp(B^4 (c))\hookrightarrow\Emb_{\om}(c,\mu)
\longrightarrow\Im\Emb_{\om}(c, \mu)
\end{equation}
where $\Symp(B^4(c))$ is the group, endowed with the $C^{\infty}$-topology, of
symplectic diffeomorphisms of the closed ball, with no restrictions on the
behavior on the boundary (thus each such map extends to a symplectic diffeomorphism of a neighborhood of $B^4(c)$ that sends $B^4(c)$ to itself), and where $\Im\Emb_{\om}(c,\mu)$ is the space of subsets of $M$ that are images of maps belonging to $\Emb_{\om}(c,\mu)$. Thus this is the set of all unparametrized balls of capacity $c$ of $M$. It is easy to show (Lemma~\ref{le:aut}) that the group $\Symp(B^4 (c))$ retracts onto its compact subgroup $\U(2)$. It is therefore equivalent, in this paper, to state the results for $\Emb_{\om}(c,\mu)$ or for $\Im\Emb_{\om}(c,\mu)$.

Denote by $\iota_c$, $c<1$, the standard symplectic embedding of $B^4 (c)$ in
$M_{\mu}$. It is defined as the composition $B^4(c)\hookrightarrow D^2(\mu- \eps)\times D^2(1-\eps)\hookrightarrow S^2(\mu) \times S^2(1) = M_{\mu}$, where the parameters between parentheses represent the areas. The results of this paper are based on the fibration
\begin{equation}\label{Fibration-2}
\Symp(M_{\mu},B_c)\hookrightarrow\Symp(M_{\mu})
\longrightarrow\Im\Emb_{\om}(c, \mu)
\end{equation}
where the space in the middle is the group of all symplectic diffeomorphisms of
$M_{\mu}$, the one on the left is the subgroup of those that preserve (not
necessarily pointwise) the image $B_c$ of $\iota_c$, and the right-hand-side arrow assigns to each diffeomorphism $\phi$ the image of the composition $\phi\circ\iota_c$.

Denoting by $\tMuc$ the blow-up of $M_{\mu}$ at the ball $\iota_c$ and by  $\Symp(\tMuc)$ the group of its symplectomorphisms, we show in Section~\ref{se:setting} that $\Symp(M_{\mu},B_c)$ is homotopically equivalent to the subgroup $\Symp(\tMuc, E)\subset\Symp(\tMuc)$ of symplectomorphisms that preserve the exceptional fiber $E$ of the blown-up space. Furthermore, we show in Section~\ref{se:setting} that the latter group $\Symp(\tMuc)$ retracts to the former $\Symp(\tMuc,E)$. We therefore have a homotopy fibration
\begin{equation}\label{Fibration-3}
\Symp(\tMuc)\hookrightarrow\Symp(M_{\mu})
\longrightarrow\Im\Emb_{\om}(c, \mu)
\end{equation}
Thus the computation of the rational homotopy type of $\Im\Emb_{\om}(c,\mu)$ can in principle be derived from the rational homotopy types of $\Symp(\tilde
M_{\mu,c})$ and $\Symp(M_{\mu})$, and the knowledge of the morphisms in the long exact homotopy sequence of~(\ref{Fibration-3}) that relies on a careful
analysis of the relations between the topology of these two groups.

Here are the results. Let us begin with the case $\mu=1$. Gromov~\cite{Gr} showed that $\Symp(M_{\mu=1})$ has the homotopy type of the semi-direct product $SO(3)\times SO(3)\times\Z/2\Z$. We prove the following in~Section~\ref{se:cohomology} as Corollary~\ref{Cor:EquivalenceG&T2}.
\begin{theorem} \label{thm:lambda=0} The group $\Symp(\tM_{\mu=1,c})$ has the homotopy type of the semi-direct product $T^2\times\Z/2\Z$.
\end{theorem}
The idea in the proof is to consider the subgroup $G$ of $\Symp(\tM_{\mu=1,c})$ which acts trivially on homology and let it  act on the (contractible !) set of configurations of triples of symplectic surfaces in classes $B-E, E, F-E \in H_2(\tM;\Z)$ (here $B$ and $F$ are, respectively, the classes of the first and second factors in $S^2\times S^2$) and show that the stabilizer ${\rm Stab} \subset G$ of a configuration $A_1\cup A_2\cup A_3$ is determined by the way each element $\phi\in{\rm Stab}$ restricts to two of these surfaces; actually, if $\phi$ is the identity on $A_1,A_3$, we show that its action on the normal bundles of $A_1$ and $A_3$ is determined up to homotopy by its action on the two normal fibers at $p = A_1\cap A_2$ and $q=A_3\cap A_2$ and that its action on the normal bundle of $A_2$ is homotopically equivalent to the fibration $\Z\to\R\to S^1$, where $S^1$ is the rotation of $A_2$ round the two fixed points $p,q$, and $\Z$ measures the torsion of the action in the normal bundle along a meridian joining $p$ to $q$ (Thus the last two actions cancel each other). Therefore $\phi$ is determined by its linear action on $A_1$ and $A_3$, that is to say, by $S^1\times S^1$. See Section~\ref{se:cohomology} for more details.

It is then easy to see the following.
\begin{cor}\label{cor:lambda=0}
When $\mu = 1$, the space $\Emb_{\om}(c,\mu=1)$ retracts, for any value of $c$, to the space
$$
\SFr(\mu=1)=\{(p,\beta)~~|~~p\in M\mbox{~and~}\beta
\mbox{~is a symplectic frame at~} p \}
$$
of symplectic frames of $M_{\mu=1}$. That is to say, $\Im\Emb(c,\mu=1)$ retracts to $M=S^2\times S^2$.
\end{cor}
\begin{proof} Consider the following commutative diagram
$$
\begin{array}{ccccc}
T^2 & \hookrightarrow & SO(3)\times SO(3) & \longrightarrow & S^2
\times S^2 \\
\downarrow &    &        \downarrow &   &   \downarrow  \\
\Symp^0(\tilde{M}_{\mu=1,c}) & \hookrightarrow & \Symp^0(M_{\mu=1}) &
\longrightarrow & \Im \Emb_{\om}(c,\mu=1).
\end{array}
$$
where the superscript $0$ denotes the identity component. Because the first two vertical arrows are homotopy equivalences, so is the last one.
\end{proof}

Consider now the case $1<\mu\le 2$. Abreu~\cite{Ab} computed the rational cohomology ring of $\Symp(M_{\mu})$ for these values of $\mu$.
\begin{theorem} [Abreu~\cite{Ab}]\label{thm:Abreu}
When $\mu\in (1,2]$, the rational cohomology of $\Symp(M_{\mu})$ is isomorphic as a ring to
$$\La (\alpha, \be_1, \be_2) \otimes S(\de)$$
where the first factor is the free exterior algebra on three generators of
respective degrees $1,3,3$ and where the second factor is the symmetric algebra on a single generator of degree~$4$.
\end{theorem}

Recall that $\la=\mu-[\mu]_-=\mu-1$. We prove in Section~\ref{se:cohomology} the following theorem, which extends Abreu's results to blown-up spaces.
\begin{theorem} \label{thm:cohomology}
When $\mu\in (1,2]$, the rational cohomology of $\Symp(\tMuc)$ is isomorphic, as a ring, to
$$\La (\alpha_1, \alpha_2 , \alpha_3) \otimes S(\eps)$$
where the first factor is the free exterior algebra on three generators all of
degree $1$, and the second factor is the symmetric algebra on a single generator which is of degree $4$ if $c$ belongs to the interval $(0,\la)$ and of degree
$2$ if $c$ belongs to the interval $[\la,1)$.
\end{theorem}

The proof is a variant of the argumentation in Abreu's work~\cite{Ab} and McDuff's work~\cite{MD:Stratification}. In some sense, it incorporates in a single
model both of their studies, that is, the case of the trivial $S^2$-bundle over $S^2$ and the case of the non-trivial one. There is, however, one new delicate point here concerning the construction of Liouville flows on open algebraic
subsets, which is resolved by approximating the symplectic form by a rational one and embedding blown-up spaces in projective spaces of arbitrarily high dimension, the dimension being determined by the rational approximation (see Lemma~\ref{le:Liouville} for details).

From Theorem~\ref{thm:cohomology}, one concludes at once the following.
\begin{cor} For any given value of $\mu\in (1,2)$, the restriction map
$\Emb_{\om}(c',\mu)\to\Emb_{\om}(c,\mu)$ is not a homotopy equivalence when $c<\la$ and $c'\ge\la$.
\end{cor}
\begin{proof} Indeed, if it were a homotopy equivalence, then, by the sequence (1), there would also be a homotopy equivalence between $\Im\Emb_{\om}(c',\mu)$ and  $\Im\Emb_{\om}(c,\mu)$ which would make the following diagram commutative:
$$
\begin{array}{ccccc}
\Symp(\tilde{M}_{\mu,c'}) & \hookrightarrow & \Symp(M_{\mu}) &
\stackrel{\rho'}{\longrightarrow} &
\Im \Emb(c', \mu) \\
  & & \downarrow \id & & \downarrow  \\
\Symp(\tMuc) & \hookrightarrow & \Symp(M_{\mu}) &
\stackrel{\rho}{\longrightarrow} &
\Im \Emb_{\om}(c, \mu).
\end{array}
$$
Hence the homotopic fibers of $\rho$ and $\rho'$ would coincide; that is,
the groups $\Symp(\tilde{M}_{\mu,c'})$ and $\Symp(\tMuc)$ would have the same homotopy type, which would contradict Theorem~\ref{thm:cohomology}.
\end{proof}

In Section~\ref{se:homotopy}, we compute the rational homotopy groups of $\Emb_{\om}(c,\mu)$ for all values of $c$ and $\mu\in [1,2]$. To achieve this, we first derive in Section~\ref{se:cohomology} the rational homotopy groups of
$\Symp(M_{\mu})$ and $\Symp(\tilde{M}_{\mu,c})$ from Theorems \ref{thm:Abreu} and \ref{thm:cohomology}: actually, the Cartan-Serre theorem \cite{CS} implies that the number of generators of the $k^{th}$ rational homotopy group (as module) is the same as the number of generators of degree $k$ in the rational cohomology (as a ring). Thus the preceding two theorems imply that the rational homotopy groups of  $\Symp(M_{\mu})$ are generated by elements $a,b_1,b_2, d$ of degree $1,3,3,4$ and that those of $\Symp(\tilde{M}_{\mu,c})$ are generated by elements $a_1,a_2,a_3$ of degree $1,1,1$ and by an element $e$ of degree $2$ or $4$ according to the value of $c$. We compute in the same section the Pontrjagin ring of $\Symp(\tilde{M}_{\mu,c})$. Note that this is, however, not  enough to derive the rational homotopy groups of $\Emb_{\om}(c,\la)$ since one must understand the morphisms in the homotopy long exact sequence corresponding to the fibration~(\ref{Fibration-3}) (see Section~\ref{se:homotopy}).

To prove that there is a homotopy jump in the spaces $\Emb_{\om}(c,\mu)$ only when $c$ reaches the value $\la$, we must show that the rational homotopy
groups of $\Emb_{\om}(c,\mu)$ remain unchanged for values of $c$ below or above the critical parameter.
\begin{theorem}  [Stability theorem] \label{Thm:Stabilite}
Given values  $\mu,\mu'\in (1,2]$, $c, c'\in (0,1)$, assume either that $c<\mu-1$ and $c'<\mu'-1$ or that $c\ge\mu-1$ and $c'\ge\mu'-1$.
\begin{enumerate}
\item[(1)] There is then a natural homotopy equivalence
$$r_{c',c}: \Emb_{\om}(c',\mu') \to \Emb_{\om}(c,\mu)$$
which coincides with the restriction map if $\mu=\mu'$ and $c'\ge c$;
\item[(2)] there is a natural homotopy equivalence between $\Im\Emb_{\om}(c',\mu')$ and\\ $\Im \Emb_{\om}(c,\mu)$; and
\item[(3)] there is a natural homotopy equivalence
$$\Symp(\tM_{\mu',c'}) \to \Symp(\tMuc).$$
\end{enumerate}
\end{theorem}

This crucial fact is proved in Section~\ref{se:inflation}. It turns out that this theorem is also useful for the determination of the rational homotopy type of $\Emb_{\om}(c,\mu)$ for $c < \la$. This is because the homotopy type of that space depends,  by the sequence (3), on the type of the space $\Symp(\tilde{M}_{\mu,c})$, which by the above theorem is homotopy equivalent to the subspace of $\Symp(M_{\mu})$ which fixes a point.

The proof of Theorem~\ref{Thm:Stabilite} is based on a beautiful idea due to McDuff in her recent paper~\cite{MD:Inflation}, which consists in considering the Lalonde-McDuff inflation technique from a global homotopical point of view.
The rest of the argument relies on the computation of Gromov's invariants in
blown-up spaces and on the analysis of the types of degeneracies of the embedded $J$-holomorphic curves associated with non-vanishing Gromov's invariants. The delicate point here is that this must be carried out for
{\em all} tamed almost complex structures, even the least generic ones (see
Section~\ref{se:inflation}).

Combining theorems \ref{thm:Abreu}, \ref{thm:cohomology} and \ref{Thm:Stabilite}, we can then compute the homotopy exact sequence of the fibration~(\ref{Fibration-3}) and conclude the following (see Section~\ref{se:homotopy}).
\begin{thm}\label{thm:rational-homotopy1}
Fix any value $\mu\in (1,2]$. Then
\begin{enumerate}
\item[(1)] the spaces $\Im\Emb_{\om}(c,\mu), c<\la=\mu-1$, are all homotopically equivalent to $S^2\times S^2$ (they have therefore the homotopy type of a finite CW-complex);
\item[(2)]  the rational homotopy groups of $\Im\Emb_{\om}(c,\mu), c\ge\la=\mu-1$, vanish in all dimensions except in dimensions $2$, $3$ and $4$ in which cases we have $\pi_2=\Q^2$, $\pi_3=\Q^3$ and $\pi_4=\Q$. (They do not have the homotopy type of a finite CW-complex.)
\end{enumerate}
\end{thm}

Let us interpret these results. We first describe what happens in $\Symp (\tM_{\mu > 1,c})$ as $c$ crosses the critical value $\la$. First, as we said above, the rational homotopy groups of $\Symp (\tMuc)$ are generated as graded module over $\Q$ by three generators $a_1, a_2, a_3$ of dimension $1$ which correspond actually to three Kahlerian $S^1$-actions that are the same for all values of $c$,  and by a fourth element $e$ which is of degree $4$ when $c < \la$ and of degree $2$ when $c \ge \la$. Denoting this fourth element by $e(4)$ and $e(2)$, respectively, it turns out that $e(2)$ is the Samelson product of two of the $1$-dimensional elements. When $c$ is smaller than $\la$, this Samelson product becomes null-homotopic and the contracting homotopy can be used to define, together with the third $1$-dimensional element, a higher Samelson product equal to $e(4)$ (see~\cite{Pinso} for the proof of this fact).

We show in Section~\ref{se:cohomology} that the one-dimensional space generated by $e$  (i.e., by $e(4)$ or $e(2)$ depending on the value of $c$) is the kernel of the map induced in rational homotopy by the inclusion $\Symp(\tMuc)\to\Diff(\tM)$ for $\mu>1$.

Concerning the spaces $\Im\Emb_{\om}(c,\mu)$, it is then easy to see the following.
\begin{prop}\label{prop:purely}
In Theorem~\ref{thm:rational-homotopy1}, the difference of rank in dimension $3$, as $c$ crosses the value $\la$, is due to the birth at that value of the new spherical generator $e(2)$ in dimension $2$ in the rational homology of $\Symp(\tilde{M}_{\mu,c})$.
\end{prop}
The difference of rank in dimension $4$ comes from the following symplectic phenomenon.
\begin{prop}\label{prop:subtle}
Fix $\mu \in (1, 2]$ and consider the generator $d$ of $\pi_4(\Symp(M_{\mu})) \otimes\Q$. Then, for each $c<\la$, there is a family of symplectic diffeomorphisms $\phi:S^4\to\Symp(M_{\mu})$ in the class $d$ that preserve some ball of capacity $c$. (This is what gives rise to the element $e(4)$.) However, such a representation is no longer possible for $c\ge\la$.
\end{prop}

Note that this proposition is compatible with the fact that the element $d$ is no longer representable  in $\pi_4(\Symp(M_{\mu}))\otimes\Q$ when $\la$ vanishes. That is to say, as the symplectic form $\mu\om_0\oplus\om_0$ on $S^2\times S^2$ approaches the form with equal weights on each factor, the class $d=[\phi]$ can be made to preserve an ``optimal ball'' of small and smaller capacity until the value $\mu =1$ is reached, where both the preserved ball and the element $[\phi]\in\pi_4(\Symp(M_{\mu>1}))\otimes\Q$ disappear.

Proposition~\ref{prop:subtle} can be considered as a nontrivial generalization of the nonsqueezing theorem in that it assigns upper bounds to higher-order widths in $M$. Let us make this more precise.
\begin{definition} Let $(M, \om)$ be a symplectic manifold of any dimension. Given a class $v \in \pi_k(\Symp (M))$, $k \ge 0$, the $v$-{\em width} of $M$ is, by definition, the supremum of the positive real numbers $c$ such that there is a symplectically embedded ball $B_c \subset M$ of capacity $c$ in M and a continuous family
$$\phi: S^k \to \Symp (M)$$
with $[\phi] = v$ and $\phi_t(B_c) = B_c$ for all $t\in S^k$. If no such ball exists, it is then easy to see that no representative of $v$ fixes even a point; in this case, we define the $v$-capacity to be zero.
\end{definition}

When $k=0$ and $v$ is the identity component of $\pi_0(\Symp (M))$, this is just the usual Gromov width, that is, the largest capacity of a ball symplectically embeddable in $M$. Proposition~\ref{prop:subtle} states that the $d$-width of $M_{\mu}$ is equal to $\mu - 1$. One may of course consider higher-order widths depending on other subspaces in $M$, though we do not pursue that question further in this paper.

We postpone to the paper~\cite{Pinso} various generalizations or extensions of
the results presented here, for instance, the treatment of the general case $\mu\in (1,\infty)$ and the proof that there is no torsion in the homotopy groups of $\Symp(\tMuc)$.

\MS
\NI
{\bf Acknowledgments.} \, F. Lalonde is grateful to Yasha Eliashberg, Silvia Anjos and Miguel Abreu for useful conversations, and M. Pinsonnault thanks Dusa McDuff for many enlightening discussions on her work with Miguel Abreu. We are also grateful to the referees for their careful reading of the paper and their pertinent suggestions. Part of this work was done while Lalonde was at Ecole Polytechnique, Palaiseau.

\section{A general framework for the study of
the space $\Emb_{\om}(B^4(c),M^4)$}\label{se:setting}

We develop in this section a general setting for the study of the space of
symplectic embeddings of the standard ball of capacity $c$ in $4$-dimensional
symplectic manifolds.

Let $\iota:B^4(c)\to M$ be any fixed symplectic embedding, let $B_c\subset M$ be its image, and let $\Emb_{\om}^{\iota}(B^4(c),M)$ be the space, endowed with the $C^{\infty}$-topology, of all symplectic embeddings of $B^4(c)$ in $(M,\om)$ that are symplectically isotopic to $\iota$ -- this is the connected component of $\iota$ inside $\Emb_{\om}(B^4(c),M)$. Let $\Im\Emb_{\om}^{\iota}(B^4(c),M)$ be the set of all images of maps in $\Emb_{\om}^{\iota}(B^4(c),M)$, with the $C^{\infty}$ topology, that is, the finest topology such that the map $\Im:\Emb_{\om}^{\iota}(B^4(c),M)\to\Im\Emb_{\om}^{\iota}(B^4(c),M)$ is continuous.
\begin{lemma}\label{le:aut}
The sequence
$$
\Symp(B^4(c),B_c)\to\Emb_{\om}^{\iota}(B^4(c),M)\to\Im
\Emb_{\om}^{\iota}(B^4(c),M)
$$
is a fibration whose fiber $\Symp(B^4(c),B_c)\simeq\Symp(B^4(c))$ retracts onto its compact subgroup~$\U(2)$.
\end{lemma}
\begin{proof}
Let us first show the second assertion. Set $G=\Symp(B^4 (c))$. Since the elements of $G$ preserve the characteristic foliation of $\p B^4(c)$, that is, the Hopf fibration, the quotient defines an element of the group $\Symp(\CP^1)$, hence a fibration
$H\to G \to \Symp(\CP^1)$ whose base retracts to $\PSU(2)$. Denote by
$\om'$ the restriction of $\om$ to $\p B^4(c)$, denote by $\Symp(\p B^4(c))$ the group of orientation preserving diffeomorphisms of $\p B^4(c)$ which preserve $\om'$, and denote by $\Symp_{\pi}(\p B^4(c))$ the subgroup that maps each leaf of the Hopf fibration to itself. Note that $\Symp_{\pi}(\p B^4(c))$ coincides with the group $\Diff^+_{\pi} (S^3)$ of diffeomorphisms of $S^3$ which map each fiber of $\pi:S^3\to\CP^1$ to itself and preserve the orientation of the fiber. This group clearly retracts to the gauge group $C^\infty(S^2,S^1)$, where $S^1$ represents the multiplication by $e^{i \theta}$ in the corresponding fiber. Thus $\Symp_{\pi}(\p B^4(c))$ retracts to $S^1$.

Now consider the restriction map
$$
H'\hookrightarrow H\stackrel{restr}{\longrightarrow}\Symp_{\pi}(\p B^4(c)).
$$
It is easy to see that this map is a fibration. Here $H'$ is the normal subgroup of $H$ which consists of all symplectic diffeomorphisms of $B^4(c)$ whose restriction to the boundary is the identity. But this group retracts to the group of symplectomorphisms with compact support in the interior of $B^4(c)$. By a result of Gromov~\cite {Gr}, it is contractible. Thus the above fibration shows that the inclusion $S^1\subset H$ given by $\theta\mapsto e^{i \theta}\id$ is a homotopy equivalence. Hence the diagram
$$
\begin{array}{ccccc}
H & \hookrightarrow & G  & \longrightarrow & \Symp(\CP^1) \\
\uparrow & & \uparrow & & \uparrow  \\
S^1 & \hookrightarrow & \U(2) & \longrightarrow & \PSU(2)
\end{array}
$$
commutes with both the left- and right-hand-side vertical arrows, being homotopy equivalences. This concludes the proof that $\Symp(B^4(c))$ retracts onto $\U(2)$.

It remains to show that the sequence of maps
$$
\Symp(B^4(c),B_c) \to \Emb_{\om}^{\iota}(B^4(c),M) \stackrel{\pi}{\to}
\Im\Emb_{\om}^{\iota}(B^4(c),M)
$$
is a locally trivial fibration. Since the fiber is the reparametrization
group, one need only show that on a $C^{\infty}$-small neighborhood $\Uu_B$ of any element $B\in\Im\Emb_{\om}^{\iota}(B^4(c),M)$, there is a continuous section of $\pi$. Because the space of embeddings considered here is connected, the extension of Hamiltonian isotopies imply that any two such embeddings differ by a symplectic diffeomorphism, and it is therefore enough to prove the existence of the section on an arbitrarily small $C^{\infty}$-neighborhood of the given ball $B_c\subset M$.  An $\eps$-neighborhood $\Uu$ of $B_c$ consists of symplectic images of $\p B^4(c)$ that are $\eps$-close to $\p B_c$ in any fixed Riemannian metric on $M$. The usual construction based on the exponential map on the normal bundle of $\p B_c$ gives a parametrization $f: \p B_c\to\p B$ of any $B\in\Uu$. The pull-back characteristic foliation is then $C^\infty$-close to the one on $\p B_c$ and one can use once again the exponential map in the normal bundle of each leaf in order to construct a foliation-preserving map $f': \p B_c\to \p B$ which therefore induces a diffeomorphism of their symplectic quotients which pulls back the symplectic form on $(\p B)/ \sim$ to a form sufficiently close to the one on $(\p B_c)/\sim$ to apply Moser's argument. Any choice of connections on both fibrations give rise to a lift of the resulting symplectic map on the quotient spaces to a symplectic diffeomorphism $f'': \p B_c \to\p B$ close to $f'$, and thus close to $f$ too. Using once again the exponential maps, there is a canonical way to extend $f''$ to a diffeomorphism $F:B_c\to B$, and by the relative version of Moser's argument, this can be made symplectic without changing the diffeomorphism on the boundary.
\end{proof}

The proof of the following lemma is entirely similar.
\begin{lemma}\label{le:aut2}
The sequence
$$
\Symp^0(M,B_c)\hookrightarrow\Symp^0(M)\to\Im\Emb_{\om}^{\iota}(B^4(c),M)
$$
is a locally trivial fibration. Here $\Symp^0(M)$ is the identity component, $\Symp^0(M,B_c)$ is the subgroup of $\Symp^0(M)$ which preserves $B_c$, and the right-hand-side arrow assigns to $\phi$ the image of $\phi\circ\iota_c$.
\end{lemma}

The only difference is that, once the section in the proof of the last lemma is
constructed, one needs to extend in some canonical way the symplectic diffeomorphism $F: B_c \to B$ to a symplectic diffeomorphism $G$ of $M$ isotopic to the identity. But this is possible because that symplectic embedding $F$ can be constructed in such a way that it is canonically isotopic to the identity. We then obtain the desired diffeomorphism $G$ by extension of Hamiltonian isotopies.

Let us now make the assumption that the space $\Emb_{\om}(B^4(c),M)$ is connected. It is known then that symplectic blow-up forms corresponding to different symplectic embeddings in $\Emb_{\om}(B^4(c),M)$ are isotopic (see~\cite{MD:BlowUps,MD:Isotopie}). Observe that in this case the sequences of Lemmas \ref{le:aut} and \ref{le:aut2} are also locally trivial fibrations when considered with the superscripts $0$ and $\iota$ suppressed.

Denote by $\Symp^{\U(2)}(M,B_c)$ the subgroup of $\Symp(M,B_c)$ that acts near $B_c$ in a $\U(2)$-linear way. Denote by $\tM_c=(\tM, \tilde{\om})$ the symplectic blow-up of $(M, \om)$ at $\iota$ and by $C_{exc}$ the exceptional
$2$-sphere, whose homology class will be denoted $E$. Write $\Symp^{\U(2)}(\tM_c,E)$ for the blow-ups of these symplectomorphisms at $\iota_c$ which therefore act linearly in a small neighborhood of the exceptional fiber $C_{exc}$, and write $\Symp(\tM_c,E)$ for the subgroup of $\Symp(\tM_c)$ that preserves (not necessarily pointwise) the exceptional fiber $C_{exc}$.
\begin{lemma}
Consider the sequence of maps
$$
\Symp(\tM_c,E) \hookleftarrow \Symp^{\U(2)}(\tM_c,E)
\stackrel{Bl}{\longrightarrow} \Symp^{\U(2)}(M,B_c) \hookrightarrow
\Symp(M,B_c)
$$
where each map is the natural inclusion except the map $Bl$ which is a
homeomorphism. All of them are homotopy equivalences.
\end{lemma}
\begin{proof} If $g$ belongs to $\Symp(M,B_c)$, its restriction to $B_c$ retracts to an element of $\U(2)$ by Lemma~\ref{le:aut}. This retraction is given by
composing $g |_{B_c}$ with some Hamiltonian $H_{t \in [0,1]}:B_c \to \R$. Now extend that Hamiltonian isotopy to $M$ so that it vanishes outside some neighborhood of $B_c$, and compose with $g$. This proves that $\Symp(M,B_c)$ retracts to $\Symp^{\U(2)}(M,B_c)$. The proof that $\Symp(\tM_c,E)$ retracts to $\Symp^{\U(2)}(\tM_c,E)$ is similar, and the other statement is obvious.
\end{proof}

We therefore have a homotopy fibration
$$
\Symp(\tM_c,E) \to \Symp(M,\om) \to \Im \Emb_{\om}(B^4(c),M)
$$
which can be used to compute the homotopy type of $\Im\Emb_{\om}(B^4(c),M)$ and thus the one of $\Emb_{\om}(B^4(c),M)$ too.

Now let us make the last assumption that the exceptional curve $E$ in
$\tilde{M}_c=(\tilde{M},\tilde{\om})$ cannot degenerate, that is, for every $\tJ$ in the space  $\Jj(\tilde{\om})$ of $\tilde{\om}$-tamed almost complex structures, there is a unique embedded $\tJ$-holomorphic representative of $E$ and no $\tJ$-cusp-curve in class $E$.
\begin{lemma}\label{le:retracts}
With these assumptions, the group $\Symp(\tM_c)$ retracts onto its\\ subgroup~$\Symp(\tM_c,E)$.
\end{lemma}
\begin{proof}  Fix any $\tilde{\om}$-tamed structure $\tJ$ for which $C_{exc}$ is holomorphic.  Let $g\in\Symp(\tM_c)$. The structure $\tJ_0= g_*(\tJ)$ has
$g(C_{exc})$ as a unique holomorphic $E$-curve. Let $\tJ_t$ be a path joining
$\tJ_0$ to $\tJ_1 = \tJ$. Note that since $\Jj(\tilde{\om})$ is contractile, the space of such paths is contractible too. By hypothesis, there is a
unique $\tJ_t$-curve in class $E$; hence, we get a one-parameter family $S_t$ of symplectically embedded {\em unparametrized} $2$-spheres joining $S_0= g_*(C_{exc})$ to $S_1= C_{exc}$.

The point is that the space of one-parameter symplectic parametrizations
$$
h_t:  S_0  \to S_t
$$
which satisfy $h_0 = \id$ (with no constraints at $t=1$), is contractible. Indeed, two such parametrizations differ by a path $g_t \in \Symp(S_0)$, and thus this space is homeomorphic to the space of paths in $\Symp(S_0)$ that start at the identity, which is obviously contractible.

Now given any such choice of parametrization, one gets a Hamiltonian isotopy $\phi_t: g_*(C_{exc}) \to \tM$ from the the identity map on $g_*(C_{exc})$ to some symplectic diffeomorphism $g_*(C_{exc}) \to C_{exc}$. Extend this isotopy to a diffeotopy $\psi_t: \tM \to \tM$ starting at the identity, and compose with $g$. The resulting symplectic diffeotopy $\psi_t \circ g: \tM \to \tM$ joins $\psi_0 \circ g = g$ to a symplectomorphism $\psi_1 \circ g$ that maps $C_{exc}$ to itself. Define
$$
\Psi : \Symp(\tM_c) \times [0,1] \to \Symp(\tM_c)
$$
by $\Psi (g,t) = \psi_t \circ g$. It depends on the following choices:
\begin{enumerate}
\item[(1)] a path $\tJ_t$  joining $\tJ_0$ to $\tJ_1 = \tJ$,
\item[(2)] a path of parametrizations $h_t$, and
\item[(3)] an extension $\psi_t$ of the Hamiltonian isotopy $\phi_t$.
\end{enumerate}
But all three choices vary in  contractible spaces; it is therefore possible to make these three choices in such a way that $\Psi$ is continuous. Note finally that, by construction, the map $\Psi$ is the identity at time $t=0$, it maps $\Symp(\tM_c)$ to $\Symp(\tM_c,E)$ at time $t=1$, and its restriction to  $\Symp(\tM_c,E) \times \{ t \}$ has image inside $\Symp(\tM_c,E)$  for all $t$. Thus $\Psi_1: \Symp(\tM_c)  \to \Symp(\tM_c,E)$ and the inclusion $\iota: \Symp(\tM_c,E) \to \Symp(\tM_c)$ are homotopy inverses.
\end{proof}

\MS

Putting the previous statements together, we get:
\begin{thm}
Let $M$ be a symplectic $4$-manifold, and let $c$ a positive real number. Suppose that:
\begin{enumerate}
\item[(1)] the space $\Emb_{\om}(B^4(c),M)$ of symplectic embeddings of the standard closed ball of capacity $c$  is non-empty and connected, and
\item[(2)] the exceptional curve that one gets by blowing up an arbitrary element of\\ $\Emb_{\om}(B^4(c),M)$ cannot degenerate in $\tM$.
\end{enumerate}
Then there is a homotopy fibration
$$
\Symp(\tM_c) \to \Symp(M,\om) \to \Im \Emb_{\om}(B^4(c),M).
$$
\end{thm}

It turns out that the manifold $M_{\mu}=(S^2\times S^2,\om_{\mu})$  satisfies these two hypotheses for all values of $\mu\ge 1$ and $c\in (0,1)$. The first one is a consequence of a theorem of McDuff~\cite{MD:BlowUps}, and the second is a computation based on the adjunction formula, the positivity of intersections, and the area constraint. Indeed, we have the following.
\begin{prop}\label{prop:unique}
Given any $\mu\ge 1$, $c\in (0,1)$, and $\tJ \in \Jj(\tilde{\om})$, there are in  $(\tMuc, \tilde{\om})$ unique embedded $\tJ$-holomorphic curves and no cusp-curve in classes $E$ and $F-E$. Moreover, these curves are Fredholm regular.
\end{prop}

Before proving the proposition, let us recall the adjunction formula. This formula states that any non-multiply covered $J$-holomorphic curve $u:\Si_g\to M$ in a class $D$ in a $4$-dimensional symplectic manifold satisfies $g_v(u) \ge g$, where $g$ is the genus of the parametrizing Riemann surface and $g_v(u)$ is the virtual genus defined by $1+\frac{D^2 - c_1(D)}{2}$, with equality if and only if the map is an embedding. For spherical classes, this means that
$$
{\mathcal Z}(C)=c_1(C)-C^2 \le 2
$$
with equality if and only if the curve $C$ is embedded. In the sequel, we will write $C$ for the image of a parametrized non-multiply covered curve, $kC$ if it is multiply covered, and $C= n_1C_1 + \ldots + n_kC_k$ for a cusp-curve
decomposition of $C$ parametrized by a tree of $k$ $2$-spheres with the restriction of $u$ to $i^{th}$ ($1 \le i \le k$) component of the tree being a $n_i$-multiple covering map onto its image.

We leave the proof of the following lemma to the reader -- it is an easy consequence of the above adjunction formula and of the positivity of $\tilde{\om}$ on $\tJ$-holomorphic curves.
\begin{lemma}\label{le:p}
If a class $pB + qF -rE \in H_2(\tM;\Z)$ admits a non multiply-covered $\tJ$-holomorphic curve, then $p \ge 0$.
\end{lemma}
\begin{proof}[Proof of Proposition~\ref{prop:unique}] The classes $E$ and $F-E$ admit embedded holomorphic representatives for generic $\tJ$-structures, and they must be regular by Hofer-Lizan-Sikorav regularity criterion (see~\cite{HLS}. By positivity of intersection, they must be unique since their self-intersection is $-1$. We must show that they cannot degenerate to cusp-curves.

Consider the case of the class $F-E$ (The case $E$ is easier and will be left to the reader.) By Lemma~\ref{le:p}, a curve in that class can only degenerate to multiples of somewhere injective curves in classes of the form $q_iF-r_iE$. The adjunction formula yields $2q_i\le 2+r_i-r_i^2$, and each value of $r_i$ different from $-1,0,1$ leads to a contradiction with the positivity of the area. The same inequality and area constraint imply that $q_i = 1$ if $r_i=1$, and $q_i=0$ if $r_i=-1$. Thus a $(F-E)$-cusp-curve must be made of positive multiples of curves in classes $F,F-E, E$; therefore the only possibility is the trivial non-cusp case.
\end{proof}
%

\section{The stability of the homotopy type of $\Emb_{\om}(c,\mu)$ below
and above the critical value $c_{crit} = \la$ } \label{se:inflation}

In this section, we prove the stability theorem (Theorem~\ref{Thm:Stabilite})
that is needed in the computation of the rational homotopy groups of
$\Symp (\tMuc)$ and of $\Emb_{\om}(c,\mu)$ in Sections 4 and 5. We recall the statement of this theorem.
\begin{theorem*}[1.6]
Given values  $\mu,\mu'\in (1,2]$, $c, c'\in (0,1)$, assume either that $c<\mu-1$ and $c'<\mu'-1$ or that $c\ge\mu-1$ and $c'\ge\mu'-1$.
\begin{enumerate}
\item[(1)] There is then a natural homotopy equivalence
$$r_{c',c}: \Emb_{\om}(c',\mu') \to \Emb_{\om}(c,\mu)$$
which coincides with the restriction map if $\mu=\mu'$ and $c'\ge c$;
\item[(2)] there is a natural homotopy equivalence between $\Im\Emb_{\om}(c',\mu')$ and\\ $\Im \Emb_{\om}(c,\mu)$; and
\item[(3)] there is a natural homotopy equivalence
$$\Symp(\tM_{\mu',c'}) \to \Symp(\tMuc).$$
\end{enumerate}
\end{theorem*}

\MS
  The rest of this section is devoted to the proof of this theorem.

\MS
(A) Let's first consider the case $c < \mu-1=\la$ and $c' < \mu'-1=\la'$.

\smallskip
To prove the three statements of the theorem, it is sufficient to prove them
\begin{enumerate}
\item[(i)] for two pairs $(\mu,c)$ and $(\mu',c')$ with $\mu= \mu'$ and
\item[(ii)] for two pairs $(\mu,c)$ and $(\mu',c')$ with $c=c'$.
\end{enumerate}
Let us establish (i), assuming, say, that $c \le c'$.

\subsection{Proof of Theorem~\ref{Thm:Stabilite} in the case (A.i) of two pairs
$(\mu,c)$ and $(\mu,c')$ with $c \le c' < \mu-1$}

Consider the following commutative diagram
$$
\begin{array}{ccccc}
\Symp(M_{\mu},B_{c'}) & \hookrightarrow & \Symp(M_{\mu}) & \to &
\Im \Emb_{\om}(c',\mu) \\
\downarrow & & \downarrow \id & & \downarrow \phi \\
\Symp(M_{\mu},B_{c}) & \hookrightarrow & \Symp(M_{\mu}) & \to &
\Im \Emb_{\om}(c,\mu)
\end{array}
$$
\NI
where as usual $\Symp(M_{\mu},B_c)$ is the space of symplectomorphisms of $M_{\mu}$ that send the standard ball $B_c \subset M_{\mu}$ to itself, where the first vertical arrow is the restriction map and the last one is induced by the quotients. It is easy to see that $\phi$ is the map that makes the following diagram commutative up to homotopy:
$$
\begin{array}{ccccc}
\Symp(B_{c'}) & \hookrightarrow & \Emb_{\om}(c', \mu) & \to & \Im
\Emb_{\om}(c', \mu)  \\
\downarrow & & \downarrow & & \downarrow \phi \\
\Symp(B_c) & \hookrightarrow & \Emb_{\om}(c, \mu) & \to & \Im
\Emb_{\om}(c, \mu)
\end{array}
$$
where the first vertical arrow is the map obtained by composing both
homotopy equivalences to $\U(2)$ and the middle one is the restriction map. Therefore, in order to show that the middle vertical arrow is a homotopy equivalence, it is sufficient to show that $\phi$ is a homotopy equivalence. Going back to the first diagram, this boils down to showing that the restriction map $\Symp(M_{\mu},B_{c'}) \to \Symp(M_{\mu},B_{c})$ is a homotopy equivalence. Consider the following diagram
$$
\begin{array}{ccccccc}
\Symp(\tM_{\mu,c'}) & \stackrel{\alpha_{c'}}{\to} & \Symp^{\U(2)}(\tM_{\mu,c'},E) & \to &
\Symp^{\U(2)}(M_{\mu},B_{c'}) & \to & \Symp(M_{\mu},B_{c'})  \\
\downarrow \psi_{c',c} &  & \downarrow r_{c,c'} & & \downarrow & & \downarrow  \\
\Symp(\tMuc) & \stackrel{\alpha_c}{\to} & \Symp^{\U(2)}(\tMuc,E) & \to &
\Symp^{\U(2)}(M_{\mu},B_{c}) & \to & \Symp(M_{\mu},B_{c})
\end{array}
$$
where all horizontal arrows are homotopy equivalences and where all vertical arrows are given by the obvious restriction or deletion maps,  except the first one, which is defined by $\alpha_c^{-1} \circ r_{c,c'} \circ \alpha_{c'}$, that is, so that the diagram commutes. Thus it is enough to show that either $\psi_{c',c}$ or $r_{c,c'}$ is a homotopy equivalence. We will consider the latter map.

It is difficult to work directly at the level of symplectic groups to construct morphisms between them and prove that they are homotopy equivalences. For this, it is natural to consider instead the fibrations
$$
\begin{array}{ccccc}
\Symp^{\U(2)}(\tM_{\mu,c'},E) & \hookrightarrow & \Diff^{\U(2)}(\tM, E) &
\stackrel{\rho'}{\to} &
\Ss (\tM_{\mu,c'},E)  \\
& & \downarrow \id  & &  \downarrow  \\
\Symp^{\U(2)}(\tMuc,E) & \hookrightarrow & \Diff^{\U(2)}(\tM,E) &
\stackrel{\rho}{\to} & \Ss (\tMuc,E) \end{array}
$$
where $\Ss (\tMuc,E)$ is the space of symplectic structures on $\tM$
diffeomorphic to the standard form $\tilde{\om}_{\mu,c}$ which coincide with the standard form near $C_{exc}$, where $\Diff^{\U(2)}(\tM,E)$ is the subgroup of diffeomorphisms of $\tM$ which act in a $\U(2)$-manner near $C_{exc}$, and where the map $\rho$ assigns to each diffeomorphism $\phi$ the push-forward $\phi_*(\tilde{\om}_{\mu,c})$ of the standard structure. Note that $\Ss(\tMuc,E)$ coincides with the space of symplectic structures on $\tM$ diffeomorphic to the standard one through a diffeomorphism which is symplectic in a neighborhood of $C_{exc}$. It is easy to see that this set coincides with the structures on $\tM$ diffeomorphic to the standard one through a diffeomorphism which is unitary in a neighborhood of $C_{exc}$. This shows that the map $\rho$ is indeed surjective. The goal is to construct the right-hand-side arrow between symplectic structures and show that it is a homotopy equivalence that makes the diagram commutative. By the homotopy theory of fibrations, this implies that the homotopy fibers are homotopy equivalent. We then finally have to show that this homotopy equivalence is homotopic to the restriction map.

It is natural to use the Lalonde-McDuff inflation technique to achieve this.
However, any direct application of that technique would require the use a family of symplectic surfaces, which belong to well-chosen classes and which vary smoothly with the symplectic form $\tilde{\om} \in \Ss (\tMuc,E)$. Moreover, one would have to show that it induces a homotopy equivalence that commutes, up to homotopy, with the action of $\Diff^{\U(2)}(\tM,E)$. This turns out to be a very intractable problem. A beautiful idea due to McDuff in her recent paper~\cite{MD:Inflation} goes round that problem in the following way. The point is that the natural space in which inflation takes place is the space $${\mathcal X}_{\mu,c}=\{(\tilde{\om},\tJ)~~:~~\tilde{\om}\in\Ss(\tMuc,E)
\mbox{~and~}\tJ\mbox{~is tamed by~}\tilde{\om}\mbox{~and is standard near~} C_{exc}\}$$
This space has two projections
$$
\begin{array}{ccc}
{\mathcal X}_{\mu,c} & \stackrel{\pi_1}{\to} & \Ss (\tMuc,E)   \\
\downarrow  \pi_2 & &  \\
\Aa_{\mu,c}
\end{array}
$$
where $\Aa_{\mu,c}$ is the space of all almost complex structures standard near $C_{exc}$ which are tamed by at least one symplectic structure in $\Ss(\tMuc,E)$. Both are Serre fibrations with contractible fibers  -- note that one has to use the additional fact that inverse images of points in the vertical one are convex to see that it is a Serre fibration. Thus these three spaces are homotopy equivalent and one has the diagram
$$
\begin{array}{ccccc}
  & & \Ss (\tM_{\mu,c'},E) & & \\
& & \uparrow \pi_1' & &  \\
  \Diff^{\U(2)}(\tM, E)
  & \stackrel{\pi'}{\to} & {\mathcal X}_{\mu,c'} &  \stackrel{\pi_2'}{\to} &
\Aa_{\mu,c'} \\
\downarrow \id & & & &    \\
\Diff^{\U(2)}(\tM, E) & \stackrel{\pi}{\to} &
  {\mathcal X}_{\mu,c} &  \stackrel{\pi_2}{\to} & \Aa_{\mu,c} \\
& & \downarrow \pi_1 & & \\
  & & \Ss (\tMuc,E) & &
\end{array}
$$
where all maps, except the fibrations $\pi$ and $\pi'$, are homotopy equivalences and $\pi_1\circ\pi =\rho, \, \pi_1'\circ\pi'=\rho'$. Here each of
$\pi$ and $\pi'$ is the action of $\Diff^{\U(2)}(\tM, E)$ by push forward on the pair consisting of the standard K\"ahler symplectic structure and the blow-up of the split complex structure.  McDuff's idea is that inflation can be used to show that the spaces $\Aa_{\mu,c}$ and $\Aa_{\mu,c'}$ are in fact {\em the same}~\footnote{In~\cite{MD:Inflation}, McDuff applies this in the simpler context of the spaces $M_{\mu}$ and $M_{\mu'}$.}. Assuming this, the homotopy fibers of $\pi_2\circ \pi$ and of $\pi_2' \circ \pi'$ would then also coincide since these two maps are equal, which would imply at once that the homotopy fibers $\Symp(\tMuc)$ and $\Symp(\tM_{\mu,c'})$ of $\pi_1\circ\pi=\rho$ and $\pi_1'\circ\pi'=\rho'$ are homotopically equivalent.

Hence the proof of case (A.i) of Theorem~\ref{Thm:Stabilite} boils down to the following two lemmas.
\begin{lemma}\label{le:coincide} For all triples $\mu,c,c'$ with $c\le c'<\mu-1$, the spaces $\Aa_{\mu,c}$ and $\Aa_{\mu,c'}$ coincide.
\end{lemma}
\begin{lemma}\label{le:restriction} The homotopy equivalence
$$
\Symp^{\U(2)}(\tM_{\mu,c'},E) \to \Symp^{\U(2)}(\tMuc,E)
$$
obtained in this way is homotopic to the restriction map.
\end{lemma}
\begin{proof}[Proof of Lemma~\ref{le:coincide}] Let us  first show that $\Aa_{\mu,c}$ is a subset of $\Aa_{\mu,c'}$. Let $\tJ\in\Aa_{\mu,c}$. Thus there is a form $\tilde{\om}$ taming $\tJ$, diffeomorphic to the standard form, and both $\tJ$ and $\tilde{\om}$ are standard in a neighborhood $V$ of $C_{exc}$. Assume for the moment that the following lemma holds. (It is proved at the end of Section~3.)
\begin{lemma} \label{le:embedded}
There are, for each almost complex structure tamed by the standard K\"ahler form $\tilde{\om}_{\mu,c}$,  embedded pseudoholomorphic $2$-spheres in the class $B+2F-E$ and in the class $B+F$.
\end{lemma}

Hence the same statement holds for our pair $(\tilde{\om},\tJ)$. Let $C$ be such an embedded $\tJ$-holomorphic sphere in class $B+2F-E$. It is $\tilde{\om}$-symplectic and has nonnegative self-intersection. By the inflation technique (see~\cite{L} or~\cite{MD:Isotopie}), one can find a representative $\be$ of the Poincar\'e dual $\PD(B+2F-E)$ with support in an arbitrarily small neighborhood $U$ of $C$, such that, for any nonnegative real number $a$, the form $\tilde{\om}+ a\be$ is symplectic and tames $\tJ$. Since $\tilde\om$ and $\tJ$ are standard on $V$, it is easy to see that we may isotop $\tilde{\om}+a\be$ in a neighborhood $V_1\subset V$ of $C_{exc}$ so that it becomes standard in some smaller neighborhood $V_2$ and tames $\tJ$ everywhere in $\tM$. Thus, in terms of the basis $\PD(B), \PD(F), \PD(-E)$ of $H^2(\tM)$, we have a form in class
$$
[\tilde{\om}_{\mu,c}] + a \PD(B+2F-E) = (1+a) \PD(B) \oplus (1+ \la +2a)
\PD(F) \oplus
(c+a) \PD(-E).
$$
Now apply inflation along a $2$-sphere $C$ in class $B+F$. For any $b\in [0, \infty)$ we then have $\tJ$ tamed by a symplectic form in class
$$
(1+a+b) \PD(B) \oplus (1+ \la +2a + b) \PD(F) \oplus (c+a) \PD(-E).
$$
Because the curve $C$ is disjoint from $C_{exc}$, this latter form remains standard near $C_{exc}$. To go back to the ratio $1+ \la$ for the relative size of the base with respect to fiber, we must choose
$$
b = \frac{a(1-\la)}{\la}.
$$
After substitution of that value of $b$ and normalization of the size of the fiber to $1$, we see that $\tJ$ is tamed by a form in class
$$
\PD(B) \oplus (1+ \la) \PD(F) \oplus  \frac{\la(c+a)}{\la+a} \PD(-E)
$$
standard near $C_{exc}$. Since $c$ is smaller than $\la$ by hypothesis, the coefficient $\frac{\la(c+a)}{\la+a}$ is larger than $c$, and as $a$ goes to infinity, that coefficient converges to $\la$. This proves that $\tJ$ is also tamed by a form in $\Ss (\tM_{\mu,c'},E)$ if $c'$ is smaller than the
upper bound~$\la$.
\begin{remark}
Note that this is a two-step process: if one wishes to achieve the inflation by adding the Poincar\'e dual to a single $\tJ$-curve, the result would then have to live in the hyperplane $y/x = 1+ \la$ of the $3$-space $H^2(\tM)$; but this is impossible if $\la$ is irrational.
\end{remark}

The proof that for all pairs $0<c\le c'<\la$, $\Aa_{\mu,c'}$ is a subset of
$\Aa_{\mu,c}$, is simpler: it is based on inflation along $2$-spheres in
classes $B+F,F$ instead of $B+2F-E,B+F$. Explicitly, if $\tJ$ is tamed by $\tilde{\om}\in\Ss(\tM_{\mu,c'},E)$, then it is also tamed by a form in class $\tilde{\om} + a\PD(B+F) +b\PD(F)$. Choosing $b = \la a$ and normalizing, we get the class $\PD(B) + (1+\la)\PD(F) + \frac{c'}{1+a}\PD(-E)$. Since $\frac{c'}{1+a}$ converges to $0$ as $a$ goes to $\infty$, this concludes the proof of Lemma~\ref{le:coincide} modulo the proof of Lemma~\ref{le:embedded} given at the end of Section~3.
\end{proof}

\MS
\begin{proof}[Proof of Lemma~\ref{le:restriction}] Denote by
$$
f_{c,c'}: \Symp^{\U(2)}(\tM_{\mu,c'},E)  \to \Symp^{\U(2)}(\tMuc,E)
$$
the map obtained in this way as homotopy equivalence between the fibers of the diagram
$$
\begin{array}{ccccc}
\Symp^{\U(2)}(\tM_{\mu,c'},E) & \hookrightarrow & \Diff^{\U(2)}(\tM, E) &
\stackrel{\rho'}{\to} &
\Ss (\tM_{\mu,c'},E)  \\
& & \downarrow \id  & &  \downarrow  \\
\Symp^{\U(2)}(\tMuc,E) & \hookrightarrow & \Diff^{\U(2)}(\tM,E) &
\stackrel{\rho}{\to} & \Ss (\tMuc,E). \end{array}
$$
The proof of Lemma~\ref{le:coincide} shows that this homotopy equivalence and its inverse are realized inside the group $\Diff^{\U(2)}(\tM, E)$ and vary continuously as the parameter $c$ runs from $c'$ to $0$ (we consider $c'$ fixed). Thus the composition
$$
f_{c,c'}^{-1} \circ r_{c,c'}: \Symp^{\U(2)}(\tM_{\mu,c'},E)  \to \Symp^{\U(2)}(\tM_{\mu,c'},E)
$$
also varies continuously. But it is equal to the identity for $c = c'$. Therefore it is always homotopic to the identity, and we conclude that $r_{c,c'}$ is homotopic to $f_{c,c'}$ for all values of $c$ smaller or equal to $c'$.
\end{proof}

\MS

\subsection{Proof of Theorem~\ref{Thm:Stabilite} in the case (A.ii) of  two pairs $(\mu,c)$ and $(\mu',c)$ with $1+c < \mu \le \mu'$}

We now prove Theorem~\ref{Thm:Stabilite} in the case (A.ii), that is, for pairs $(\mu,c), (\mu',c)$ satisfying $1+c < \mu \le \mu'$.  As in the case (i), the proof of all three statements of the Theorem boils down to establishing the existence of a natural homotopy equivalence between
$\Symp^{\U(2)}(\tM_{\mu',c},E)$ and $ \Symp^{\U(2)}(\tMuc,E)$ or, equivalently,
between $\Symp(\tM_{\mu',c})$ and  $\Symp(\tMuc)$. Considering the diagram
$$
\begin{array}{ccccc}
\Symp(\tM_{\mu',c}) & \hookrightarrow & \Diff(\tM) &
\stackrel{\rho'}{\to} &
\Ss (\tM_{\mu',c})  \\
& & \downarrow \id  & &  \downarrow  \\
\Symp(\tMuc) & \hookrightarrow & \Diff(\tM) &
\stackrel{\rho}{\to} & \Ss (\tMuc) \end{array}
$$
it is enough to show that the space $\Aa_{\mu',c}$ is equal to $\Aa_{\mu,c}$, where this time $\Aa_{\mu,c}$ is the space of almost complex structures on $\tM$ that are tamed by at least one symplectic structure diffeomorphic to the standard symplectic form $\tilde{\om}_{\mu,c}$.  This is easy. Let $\tJ \in \Aa_{\mu,c}$. Inflation along a curve in class $F$ gives a structure tamed by a form in class $PD(B) + (\mu+a)PD(F) + c PD(-E)$. Hence, choosing $a = \mu'-\mu$, we see that $\tJ \in \Aa_{\mu',c}$. Conversely, if $\tilde \in \Aa_{\mu',c}$, inflation along a curve in class $B+F$ means that $\tJ$ is tamed by a form in class
$$
(1+a)PD(B) + (\mu'+a)PD(F) + cPD(-E)
$$
or, after normalization,
$$
PD(B) + \frac{\mu' + a}{1+a} PD(F)  + \frac{c}{1+a} PD(-E).
$$
Because $\mu\le\mu'$, we can find a value of $a$ such that the resulting form belongs to  $\Ss(\tM_{\mu,c_0})$ with $c_0\le c$. Applying case A(i), we can increase $c_0$ to $c$ keeping $\mu$ fixed. Thus $\tJ$ belongs to $\Ss(\tMuc)$.
\qed

The proof of case (B), i.e the the case $c\ge\mu-1=\la$ and $c'\ge\mu'-1=\la'$, is easier than case (A) (because there are in this case only two strata -- see Sections 3.3 and 4.1). It is left to the reader.

\subsection{Existence of pseudo-holomorphic curves in required classes}

To complete the proof of Theorem~\ref{Thm:Stabilite}, there only remains to establish Lemma~\ref{le:embedded}.

\begin{proof}[Proof of Lemma~\ref{le:embedded}] Let $\tJ$ be tamed by $\tilde{\om}_{\mu,c}$ with $0<c <\la\le 1$. We must show that there are embedded $\tJ$-holomorphic $2$-spheres in classes $B+2F-E$ and $B+F$.

Let us first establish a few general facts. The first one, concerning the Gromov invariants in blown-up spaces, is elementary.  Recall that for a $4$-dimensional symplectic manifold $X$ and a class $D \in H_2(X;\Z)$ such that $k(D)=\frac{1}{2}(c_1(D) + D^2)\ge 0$, the Gromov invariant of $D$ as defined by Taubes in~\cite{Ta} counts, for a generic $J$ tamed by $\om$, the algebraic number of embedded $J$-holomorphic curves in class $D$ passing through $k(D)$ generic points.
\begin{prop}\label{prop:SW}
Let $(X,\om)$ be a closed $4$-dimensional symplectic manifold, and assume that some class $D$ has non-vanishing Gromov invariant. Let $(\tilde X, \tilde{\om})$ be the symplectic blow-up of any symplectic embedding of a ball. Then the class $D$ also has non-vanishing Gromov invariant in $(\tilde X,\tilde{\om})$. The same is true of the class $D-E$ if the necessary condition $k(D-E) \ge 0$ is satisfied.
\end{prop}
\begin{proof} First note that the Gromov invariant depends only on the deformation class of the symplectic form. Therefore one may assume that the blow-up is made at some ball of arbitrarily small radius. Now consider first the case of class $D$. We may suppose that the embedded $J$-curves $C_1,\ldots,C_{\ell}$ in class $D$ passing through the $k(D)$ points $p_1,\ldots, p_k$ are disjoint from the center $p$ of the ball that we assume to be a generic point. Since $J$ is
generic, each of these curves is $J$-regular. They are therefore regular for some almost complex structure $J'$ which coincides with $J$ everywhere except in a small neighborhood of $p$ where $J'$ is integrable. Moreover, there must be such a structure $J'$ for which the only $J'$-curves passing through these $k$ points are the above $C_i$'s. Indeed, if not, there would be a sequence of such structures $J'_n$ coinciding with $J$ on $X - B_p(1/n)$ and having a
holomorphic curve $A_n$ passing through the above $k$ points and some point $q_n\in B_p(1/n)$. By Gromov's compactness theorem, this would give a
$J$-cusp-curve in class $D$ passing through the $k+1$ generic points $p_1,\ldots, p_k, q_n$. But this is impossible since the space of such $D$-cusp-curves is smaller than the dimension of embedded $D$-curves except if the manifold admits pseudoholomorphic curves with negative Chern class with respect to generic almost complex structures. But this cannot happen in dimension $4$.

Now blow-up the locally integrable structure $J'$ to get a structure
$\tJ'$ on $\tilde X$. Its only representatives in class $D$ passing through
$p_1,\ldots, p_k$ are the proper transforms of the $C_i$'s, which are
all regular. Indeed, if there were another one, it would not be a cusp-curve
for dimensional reasons and would not meet the exceptional fiber either (because $D\cdot E = 0$), therefore its blow-down in $X$ would give another $D$-curve, a contradiction. Thus the Gromov invariant of $D$, as a class inside $\tilde X$, must be equal to the Gromov invariant of $D$ seen inside $X$.

Consider now the case $D-E$. By hypothesis, $k(D-E)\ge 0$, which implies that $k(D)\ge 1$. There are, say, $\ell$ embedded regular $J$-curves of $X$ in class
$D$ passing through generic points $p_1,\ldots, p_{k-1},p$. By regularity,
these curves persist and remain regular for a nearby almost complex structure $J'$ integrable near $p$. We claim that there is no other $J'$-curves in class $D$ passing through the same points. Assuming the contrary, by the compactness theorem, we get a sequence of structures $J_n'$ converging to $J$ with at least two $J_n'$-curves converging to the same regular $J$-curve. But for regular curves, such a bifurcation is impossible. In the blow-up space with the blow-up almost complex structure, these curves lift to regular embedded curves in class $D-E$ passing through $k(D)-1 = k(D-E)$ generic points of $\tilde X$. There cannot be any other such curve because it would have to meet the exceptional divisor transversally at one point and therefore its blow-down would lead to a contradiction. Thus, once again, the Gromov invariant is the same.
\end{proof}

We note, by the way, that this proposition has the following useful corollary.
\begin{cor}\label{cor:capacite} Let $(X, \om)$ be a closed $4$-dimensional symplectic manifold, and denote by ${\rm cap} (X, \om)$ the supremum of the capacities of balls symplectically embeddable in $M$. Let $D$ be any integral homology $2$-class with $k(D)\ge 1$ and area less than ${\rm cap} (M,\om)$. Then the Gromov invariant of $D$ vanishes.
\end{cor}
\begin{proof} If not, the previous proposition would imply that the Gromov invariant of $D-E$ would not vanish in $(\tilde X,\tilde{\om}_c)$, the blow-up of $(X,\om)$ at some small ball of capacity $c$. But as $c$ increases, the Gromov invariant remains unchanged because it only depends on the deformation class of the symplectic structure. For large values of $c$, the symplectic area of $D-E$ becomes negative, a contradiction.
\end{proof}

\begin{remark} Of course, the contraposition of Corollary~\ref{cor:capacite} is equally useful: if one knows that $D$ has non-vanishing Gromov invariant and satisfy $k(D) \ge 1$, then ${\rm cap} (X, \om)$ is bounded from above by the area of $D$. This is simply a mild generalization of the proof of the nonsqueezing theorem since one may replace the $S^2$-factor by any class with nonvanishing Gromov invariant that has at least one $J$-representative through each point.
\end{remark}

\noindent Let us come back to the proof of Lemma~\ref{le:embedded}. We need the following facts.

Fact 1. Recall that Lemma~\ref{le:p} states that any class $pB+qF-rE$  that admits nonmultiply covered $\tJ$-holomorphic representatives must satisfy $p \ge 0$. Thus if a curve in class $B+qF-rE$ degenerates, the resulting cusp-curve must be the union of a curve of the form $B+\ldots$ with curves of the form $qF-rE$.

%
Fact 2. Because $E$ is always represented (see Proposition~\ref{prop:unique}), the positivity of intersection implies that the coefficient $r$ is always non-negative, except of course in the case $p=q=0$.

%
Fact 3. The adjunction formula for genus zero nonmultiply covered curves implies that a curve of the form $B+qF-rE$ satisfies $r=0$ or $1$ and must be embedded, and a non-multiply-covered curve of the form $qF-rE$ satisfies $2q-r(1-r)\le 2$ and is embedded if $q$ is nonnegative.

We can now prove Lemma~\ref{le:embedded}. By the wall-crossing formula for the Sieberg-Witten invariant in Li and Liu~\cite{LL}, the Gromov invariant of all classes $D$ of $M=S^2\times S^2$ with $k(D)\ge 0$ is equal to $1$. Thus by Proposition~\ref{prop:SW}, the Gromov invariants of both $B+F$ and $B+2F-E$ in $\tM$ do not vanish. Hence there are embedded $J$-curves in class $D=B+2F-E$ or $D=B+F$ for generic $\tJ$'s passing through any set of $k(D)$ generic points. Let $\tJ$ be any almost complex structure tamed by $\tilde{\om}_{\mu,c}$. By the compactness theorem, there are (embedded or not) $\tJ$-curves in class $D$ passing through any set of $k(D)$ points; in order to show that there are embedded $\tJ$-curves in class $D$, it is enough to show that the dimension of unparametrized cusp-curves is less than the dimension of the moduli space corresponding to the conditions given by $k(D)$ points. Consider first the case $B+2F-E$. First assume that $\tJ$ belongs to the stratum $\tilde{\Jj}_0$ admitting curves in classes $B$; then by positivity of intersection, no component can realize a class $pB + qF -rE$ with $q$ strictly negative. Hence all components have coefficient $q=0,1$ or $2$. Combining this with Facts 1, 2, and 3, this implies that an embedded $B+2F-E$ curve can only degenerate to a cusp-curve that is the union of an embedded curve in class $B+qF-rE$ with $q\in\{0,1,2\}$, $r\in\{0,1\}$, and of positive multiples of embedded curves in classes of the form $qF-rE$, with $(q,r)=(0,-1)$, $(1,0)$, or $(1,1)$. In all cases, each curve $C_i$ in the decomposition is embedded with strictly positive first Chern class.  Hence  by the regularity criterion of Hofer, Lizan, and Sikorav~\cite{HLS}, the dimension of the moduli space of each component is equal to the one predicted by the index formula. Since none of the $n_i C_i$ in the decomposition is a multiple of a class of negative Chern number, we conclude by a simple counting argument that the dimension of each cusp-curve type is strictly less than the dimension of the moduli space of the embedded curves in class $B+2F-E$. Hence, there must be a $\tJ$-embedded curve in that class.

If $\tJ$ belongs to the stratum $\tilde{\Jj}_1$ of structures admitting $(B-F)$-curves, then the positivity of intersection of $B-F$ with $B+qF-rE$ and with $qF-rE$ implies once again that $q \ge 0$, and the rest of the argument is the same. Finally, if $\tJ$ belongs to the stratum $\tilde{\Jj}_2$ of structures admitting $(B-F-E)$-curves, then the positivity of intersection of $B+qF-rE$ with $B-F-E$ implies that $q$ is positive (except in the case $B+qF-rE=B-F-E$), and the positivity of intersection of $qF-rE$ with both $E$ and $B-F-E$ implies that $q$ is non-negative. We can then proceed as before.
\end{proof}

\noindent Thus Theorem~\ref{Thm:Stabilite} is proved.\qed

\section{The cohomology ring and the Pontrjagin ring
of~$\Symp(\tMuc)$} \label{se:cohomology}

\subsection{Stratification of the space of almost complex
 structures and stabilizers}

In this section, we compute the rational cohomology ring of $\Symp(\tMuc)$ using techniques similar to the ones developed by Abreu in \cite{Ab}, Abreu and McDuff in \cite{AM}, and McDuff \cite{MD:Stratification}. Because $\tMuc$ can be blown down to either $S^2 \times S^2$ or $\CP^2 \#\bar{\CP^2}$ (by blowing down either the class $E$ or the class $F-E$), the study of the group $\Symp(\tMuc)$ incorporates in some sense the two different studies in Abreu and McDuff~\cite{AM} on the symplectic groups of  $S^2 \times S^2$ and $\CP^2 \#\bar{\CP^2}$. We explain in this section how the two geometries fit together, and refer the reader to Abreu and McDuff in all cases when the proofs are mild adaptations of those contained in~\cite{Ab}, \cite{AM}, and~\cite{MD:Stratification}.

We first explain the general framework for the calculation of that cohomology ring of $\Symp(\tMuc)$ for all values $\mu \ge 1$, and then restrict ourselves later to the case $\mu \le 2$.

\begin{remark} \label{re:rationality} By the stability theorem (Theorem~\ref{Thm:Stabilite}), the homotopy type of $\Symp(\tMuc)$ does not depend on $\mu, c$ as long as $c$ remains below or over the critical parameter $\la = \mu-1$. We may therefore assume without loss of generality that $\mu$ is rational; this is convenient since the proofs of some of the basic propositions of this section are most naturally carried out under this assumption. We always state it explicitly when the hypothesis ``$\mu$ rational'' is needed.
\end{remark}

For $k \ge 0$, denote by $D_{2k} \in H^2(\tMuc,\Z)$ the class $B-kF-E$ and by $D_{2k+1}$ the class $B-(k+1)F$.  We already know from the Section~\ref{se:inflation} that the classes $E$ and $F-E$ cannot degenerate. However the class $B-E$ can and does actually degenerate when $\mu$ is larger than $1$. Recall that $[\mu]_-$ is the lower integral part of $\mu$, that is to say, the largest integer amongst all those that are strictly smaller than $\mu$, and $\la = \mu - [\mu]_-$.
\begin{prop}\label{Prop:StructureCourbesCusp}
Let $\mu \ge 1$ and $0 < c < 1$. For each $\tJ\in\tjj(\omlc)$, there is one and only one of the classes $\{D_i\}_{0 \le i \le \ell}$ that is represented by a $\tJ$-holomorphic curve.  Moreover that curve is embedded and unique in its homology class. Here $\ell$ is equal to $2[\mu]_-$ if $c < \la$ and is equal to $2[\mu]_- - 1$ if $c\ge\la$.
\end{prop}
\begin{proof}
This is an easy consequence of the analysis of degeneracies of the unique embedded $\tJ'$-curve in class $B-E$ for a generic $\tJ'$, as $\tJ'$ approaches any given  almost complex structure $\tJ$ compatible with $\tilde{\om}_{\mu,c}$.

By Gromov's compactness theorem, the exceptional $\tJ'$-curve $C$ in class $B-E$ must degenerate to a $\tJ$-cusp-curve as $\tJ'$ approaches $\tJ$. Recall that Lemma~\ref{le:p} states that any class $pB + qF -rE$  which admits nonmultiply covered $\tJ$-holomorphic representatives must satisfy $p \ge 0$. Thus if a curve in class $B-E$ degenerates, the resulting cusp-curve must be the union of a curve of the form $B+\ldots$ with curves of the form $qF - rE$.

Because $E$ is always represented (Proposition~\ref{prop:unique}), the positivity of intersection implies that the coefficient $r$ is always nonnegative, except of course in the case $p=q=0$. The adjunction formula for genus zero nonmultiply covered curves implies that a curve of the form $B + qF - rE$ satisfies $r=0$ or $1$ and must be embedded. But if the coefficient $q$ were strictly positive, there would be some other component of the form $-q'F -r'E$ with $q'> 0$ and $r'\ge 0$, but this is impossible since this would have negative area. The fact that only one such class is represented for each $\tJ$ (as well as the fact that the $\tJ$-curve is unique in its homology class) is due to the positivity of intersection; it is embedded by the adjunction formula. For more details, see the proof of Lemma~\ref{le:embedded} in Section~\ref{se:inflation}, where a similar but more difficult statement is proved.
\end{proof}

Corresponding to each degeneracy of type $D_i$, there are
\begin{enumerate}
\item a subset $\tjj_i \subset \tjj(\omlc)$ consisting of all almost complex
structures $\tJ$ such that there is a $\tJ$-curve in class $D_i$;
\item a configuration space $\Cc_i$ consisting of all triplets $(A_1,A_2,A_3)$ formed of $\tilde{\om}_{\mu,c}$-symplectic embedded surfaces in $\tMuc$ which intersect transversally and positively, and lie in classes $F-E, E, B-kF-E$, respectively, if $i = 2k$ and in classes $E, F-E, B-(k+1)F$, respectively, if $i=2k+1$ (note that, in all cases, $A_2$ intersects each of $A_1$ and $A_3$ at one point, and $A_1 \cap A_3 = \emptyset$);
\item a K\"ahler structure $\tJ_i \in \tjj_i$ for which the group of $(\tilde{\om}_{\mu,c}, \tJ_i)$-K\"ahler isometries is a $2$-torus (see below, after the explanation of Proposition~\ref{Prop:Stratification} ).
\end{enumerate}

As in~\cite{AM}, we will denote by $\tjj_{0 1 \ldots i}$ the union $\tjj_0 \cup \ldots \tjj_i$. Proposition~\ref{Prop:StructureCourbesCusp} shows that the set $\{\tjj_i\}_{0 \le i \le\ell}$ forms a partition of $\tjj(\omlc)$. The index formula and standard arguments from the theory of $J$-curves imply that $\tjj_0$ is an open dense subset of $\tjj(\omlc)$ and that each stratum $\tjj_i$ has  finite codimension equal to $2i$ in $\tjj(\omlc)$. It is, however, more delicate to see that these subsets fit together to form a stratification in the sense of the following proposition. This was proved by McDuff in \cite{MD:Stratification} in a slightly different context. However her arguments carry over to our case with obvious changes. This gives the following.
\begin{prop}\label{Prop:Stratification}
The space $\tjj(\omlc)$ is equal to the disjoint union $\tjj_0\cup\cdots\cup\tjj_\ell$, where each $\tjj_i$ is a nonempty, co-orientable, and smooth submanifold of codimension $2i$. The closure of $\tjj_i$ in $\tjj(\omlc)$ is the union $\tjj_i\cup\cdots\cup \tjj_\ell$. Moreover, each $\tjj_i$ has a neighborhood $\Nn_i \to \tjj_i$ in $\tjj$ which, once given the induced stratification, has the structure of a locally trivial fiber bundle whose typical fiber is a cone over a finite dimensional stratified space called the {\em link} of $\tjj_i$ in $\tjj$. Finally, the link $\Ll_{i}$ of each $\tjj_{i}$ in $\tjj_{i-1}$ is a circle.
\end{prop}

\BS

Here is the idea of the proof that shows why the link is a circle. As explained above, to each almost complex structure $\tJ$ in the stratum $\tjj_{i}$, one assigns the unique configuration of $\tJ$-holomorphic curves $A_1\cup A_2\cup A_3$. It is actually shown below (Proposition~\ref{Prop:EquivalenceJi&Ci}) that the correspondence between strata of almost complex structures and configuration spaces is a homotopy equivalence. Denote by $p$ the point where the curves $A_2$ and $A_3$ intersect positively and transversally. Let's consider the gluing of the $\tJ$-curves  $A_2$ and $A_3$. Consider this first from a homological point of view only. If $i=2k$, this gives a curve in class $[A_2 + A_3] = (B-kF-E)+E = B-kF$, which is  the class $[A_3] = D_{2k-1}$ of a configuration corresponding to the stratum $\tjj_{2k-1} = \tjj_{i-1}$. If $i=2k + 1$, this gives a curve in class $[A_2 + A_3] = (B-(k+1)F)+(F-E) = B-kF-E$, which is once again the class $[A_3] = D_{2k}$ of the configuration corresponding to the stratum $\tjj_{2k} = \tjj_{i-1}$.

Using the gluing theorem of transversal rational curves, one can resolve the singular node $p$ in order to get a curve $C$ in the class $[A_2 + A_3]$. The point is that this new curve $C$ is holomorphic with respect to some almost complex structure $\tJ'$ in the stratum $\tjj_{i-1}$. In other words, this gluing process can be seen as a partial converse to the Gromov compactness theorem since it describes a neighborhood of $A_2\cup A_3$ in the compactified moduli space $\overline\Mm([A_2]+[A_3],\tjj_{i-1})$. This gluing process depends continuously on $\tJ\in\tjj_{i}$ and on the choice of some vector $a$ in a neighborhood $V_{\tJ}$ of the zero element in the 1-dimensional complex vector space $T_p(A_2)\otimes T_p(A_3)$. The analysis done in~\cite{MD:Stratification} then shows (i) that one may use this gluing to identify $V_{\tJ}$ with a normal neighborhood of $\tJ$ in $\tjj_{i-1}$, and (ii) that moreover,  these identifications can be made coherently so that they yield the structure of a topological disc bundle filling a closed neighborhood of $\tjj_{i}$ in $\tjj_{i-1}$.

\MS \MS
The goal of the next propositions is to establish the relations between (1), (2) and (3) above, namely, that
\begin{itemize}
\item[--] each stratum $\tjj_i$ contains a K\"ahler structure $\tJ_i$ obtained from the blow-up of the Hirzebruch surface $W_i$ (this applies for both odd and even Hirzebruch surfaces) and that as $i$ varies from $0$ to $\ell$, the K\"ahler structures run over all toric geometries on $\tM_{\la,c}$;
\item[--] each stratum is homotopically equivalent to the configuration space $\Cc_i$, and although $\Symp(\tilde{M}_{\mu,c})$ does not act transitively on $\tjj_i$ (since, for instance, the action preserves the  (non) integrability of $\tJ$ or any Riemannian invariant of the pair $(\om, \tJ)$) or on $\Cc_i$ (there are moduli at the intersection points), it does act transitively on the subspace $\Cc_i^0$ of orthogonal configurations;
\item[--] the stabilizer of this action can be identified with the stabilizer of the action of the group on $\tjj_i$, that is to say to the $2$-torus $\tT^2_i$ of K\"ahler isometries of $\tJ_i$. Moreover, its rational homology injects in the rational homology of $\Symp(\tilde{M}_{\mu,c})$.
\end{itemize}

\BS

Thus $\Symp(\tilde{M}_{\mu,c})$ contains all the tori $\tT^2_i, 0 \le i \le \ell$, and it is not difficult to identify their pairwise intersections in $\Symp(\tilde{M}_{\mu,c})$. As we will see below, this implies that the rational cohomology of $\Symp(\tilde{M}_{\mu,c})$ is the tensor product $H^*(\tT^2_i)\otimes H^*(\tjj_i)$ for each $0 \le i \le \ell$. Since the strata must fit together to form the contractible space $\tjj$, these $\ell + 1$ equations, together with the Leray theorem on the ring structure of $H^*(\Symp(\tilde{M}_{\mu,c}))$,  lead to an effective way of computing $H^*(\Symp(\tilde{M}_{\mu,c}))$.

\MS
We first briefly recall the definition of the Hirzebruch surfaces.  For any $\mu > 0$ and any integer $i \ge 0$ satisfying $\mu - \frac i2 > 0$, let $\CP^1 \times\CP^2$ be endowed with the K\"ahler form $(\mu - \frac i2) \tau_1 + \tau_2$ where $\tau_{\ell}$ is the Fubini-Study form on $\CP^{\ell}$ normalized so that the area of the linear $\CP^1$'s is equal to $1$. Let $W_i$ be the corresponding Hirzebruch surface, that is, the K\"ahler surface defined by
$$
W_i=\{([z_0,z_1],[w_0,w_1,w_2])\in\CP^1\times\CP^2 ~|~ z_0^iw_1=z_1^iw_0 \}
$$
It is well-known that the restriction of the projection $\pi_1:\CP^1\times\CP^2\to\CP^1$ to $W_i$ endows $W_i$ with the structure of a K\"ahler $\CP^1$-bundle over $\CP^1$ which is topologically the trivial one $S^2 \times S^2$ if $i$ is even and is the  nontrivial one $S^2\times_{\tau} S^2=\CP^2\#\bar{\CP}^2$ if $i$ is odd. In this correspondence, the fibers are preserved and the zero section of this bundle
$$
s_0 = \{([z_0,z_1],[0,0,1])\}
$$
corresponds to the section  of self-intersection $-i$ that lives in $S^2\times S^2$ if $i$ is even and in the nontrivial $S^2$-bundle if $i$ is odd. Thus it represents the class $B-(i/2)F$ (resp. $\si_{-1} -(\frac{i-1}{2})F$, where $\si_{-1}$ is the section of $S^2 \times_{\tau} S^2$ of self-intersection $-1$). By the classification theorem of ruled symplectic $4$-manifolds, this correspondence establishes a symplectomorphism between $W_i$ and  $M_{\mu}$ for all even $i$'s and between $W_i$ and $M_{\mu}'$ for all odd $i$'s. Here $M_{\mu}$ is $S^2\times S^2$ endowed with the form $\om_{\mu}$ and $M_{\mu}'$ is $S^2 \times_{\tau} S^2$ with the form that gives area $1$ to the fiber and area $\mu$  to the ``base'' $(\si_{-1} +\si_1)/2$. In the case of $S^2 \times S^2$ with the symplectic form $\om_{\mu}$, this means that we have, for each small enough even number $i$ (i.e $i/2 < \mu$), a K\"ahler structure on $M_{\mu}$, coming from $W_i$, having the section $B-(i/2)F$ holomorphically represented. A similar comment applies to $M_{\mu}'$.

We refer to Audin~\cite{Au} for a description of the groups $K_i$ of isometries of these K\"ahler structures. For the manifold $W_{0}$, the group $K_0$ is the product of $SO(3)\times SO(3)$, while in all other even cases $M_{2k}, k > 0$, it is  $SO(3) \times S^1$. Here the $S^1$-factor is the rotation in the fibers of $W_{2k}$ by complex multiplication (and it therefore fixes pointwise the zero section), while the $SO(3)$-factor is a lift to $W_{2k}$ of the standard $SO(3)$-action on the base. For the odd $W_{2k+1}$, the K\"ahler group $K_i$ is $\U(2)$ seen as the quotient $SU(2) \times S^1 / \{(1, \pm 1), (-1, \pm \sqrt{-1}) \}$ (see~\cite{Au} or~\cite{AM} for a precise description of this action). We will only need the following easy consequence of the description in~\cite{Au}.
\begin{prop} \label{prop:rotations} Given a fiber $F_i$ in $W_i$, the subgroup of $K_i$ which preserves the configuration $s_0\cup F_i$ is equal to a 2-torus $T^2_i$, where the second $S^1$-factor $t_i$ acts by rotations in the fibers and the first one $s_i$ is a fiberwise diffeomorphism obtained in the even case as the subgroup of $SO(3)$ that preserves $F_i$, and in the odd case as the quotient of the diagonal elements $e^{i\theta}\id$ of the $SU(2)$-factor.
\end{prop}
\MS
Because $s_0$ is an exceptional curve of $W_i$, it is preserved by any K\"ahler isometry. On the other hand, since all these K\"ahler isometries are fiberwise diffeomorphims, preserving a given fiber amounts to preserving the intersection point with $s_0$. Thus this proposition states that the subgroup of $K_i$ that preserves a point on $s_0$ is $T^2$. We leave to the reader to check that this easily implies the following.
\begin{prop} The group of K\"ahler isometries of the blow-up $\tilde W_{i,c}$ of $W_i$ at a standard ball of any capacity $c<1$ centered at a point on $s_0$ can be identified with the 2-torus $T^2_i$ of Proposition~\ref{prop:rotations}.
\end{prop}
\BS
Note that the blow-up $\tilde W_{i,c}$ is symplectically diffeomorphic to $\tMuc$ {\em regardless} of the parity of $i$. In the even case, the configuration $F_i, p, s_0$ is sent, by the blow-up operation centered at the point $p = F_i \cap s_0$, to the configuration $F-E, E, D_i$, while in the odd case the same configuration is sent to $E, F-E, D_i$. Note that the $F_i$'s live in $W_i$ but the classes $E,F,D_i$ refer to the standard classes of $\tMuc$ as defined in the Introduction. Therefore, each torus $T^2_i$ gives rise to an abelian subgroup of $\Symp(\tMuc)$ which we denote by $\tT^2_i$. Let us now express this in terms of the precise description of the Hirzebruch surfaces given above. Consider the following $T^2_i$-action on $W_i$:
\begin{eqnarray}\label{Action1}
(\theta_1,\theta_2)\cdot([z_0,z_1],[w_0,w_1,w_2]) & = & ([\theta_1z_0,z_1],[\theta_1^{i}w_0,w_1,\theta_1^{i}\theta_2w_2])
\end{eqnarray}
where the action of $\theta_j$ is the product by $\exp{2 \pi \sqrt{-1} \theta_j}$. Identify the fiber $F_i$ with the fiber of $W_i$ over the point $[z_0,z_1]=[1,0]$. The above action then coincides with the action of Proposition~\ref{prop:rotations}.  Blow up this action at the fixed point  $s_0 \cap F_i$ (note that as a rule, except when stated otherwise, all blow-ups in the sequel are taken at that point  $s_0 \cap F_i$ which, in the case of $W_0$, is identified with the center of the standard ball $\iota_c$ in $M_{\mu}$). As in Proposition~\ref{prop:rotations}, denote by $s_i, t_i$ the generators of the $T^2_i$-action corresponding to $\theta_1, \theta_2$, respectively, and by  $x_i,y_i$ the generators of the fundamental group of $\tT^2_i$ which are the blow-ups of $s_i, t_i$. It is then straithforward to check the following.

\begin{prop}\label{prop:action}
The first generator $x_i$ of $\tT^2_i$ (corresponding to $\theta_1$) rotates the surface $A_3$ of the configuration $C_i$ round the point $q = A_2 \cap A_3$ but fixes pointwise the surface $A_1$, and the second generator $y_i$ (corresponding to $\theta_2$) rotates the surface $A_1$ of the configuration $C_i$ round the point $p = A_1 \cap A_2$ and fixes $A_3$ pointwise. Because this action comes from one that fixes a point corresponding to the blow-down of $A_2$, the restriction of the action to $A_1 \cup A_2$ (or to $A_2 \cup A_3$) rotates the two surfaces independently.
\end{prop}

 \MS
  
The toric action $T^2_i$ on $W_i$ can be blown up at any of its fixed points. Thus one gets another torus action by blowing up $T^2_i$ at the intersection of $F_i$ with the section at infinity $s_{\infty}$. Let's denote by $(\tT^2_i)_{\infty}$ this action and by $(\tT^2_i)_{0} = \tT^2_i$ the former one.  

\begin{prop} \label{prop:invert} We have the relations
$$
(\tT^2_{2k-1})_{0} = (\tT^2_{2k})_{\infty}  
$$
and
$$
 (\tT^2_{2k-1})_{\infty}= (\tT^2_{2k-2})_{0}
$$
among the toric actions on $\tMuc$ which hold up to the following reparametrisation of the tori: in both cases, the $y_i$'s generators correspond, and the generator $x_{2k-1}$ is mapped to $x_{2k} - y_{2k}$ in the first case and to $x_{2k-2} + y_{2k-2}$ in the second case.
\end{prop}

\begin{proof} In both cases, the birational correspondence $\phi$ between $W_{2k-1}$ on the one hand and $W_{2k}$ or $W_{2k-2}$ on the other hand is given by blowing up $W_i$ at the point $*$ (equal to $s_0 \cap F_i$ in the first case and to $s_{\infty} \cap F_i$ in the second case) and then blowing down along the proper transform of the fiber passing through the point $*$. This birational transformation between odd and even $W_i$'s obviously preserves the fibers:  it is a lift to the Hirzebruch surfaces of the identity map on the base $\CP^1$ of all $W_i$'s. Thus in both equalities, the correspondence preserves the $t_i$'s. However, in the first case, the $s_{2k-1}$ circle action fixes the fiber over $[1,0]$ pointwise and rotates the fiber positively over $[0,1]$ by $2k-1$ turns. Because the transform $\phi$ does not affect the latter fiber, the same is true of the image action $\phi_{*}(s_{2k-1})$ on that fiber. But this image rotates the fiber over $[1,0]$ by minus one turn. Note that the difference of full one-turn rotations, at the target space $W_{2k}$, between the two fibers is equal to $2k$ as it should be since, as one sees in formula (\ref{Action1}), every circle action inside $T^2_i$ whose projection onto the base $\CP^1$ agrees with the projection of $s_i$ (i.e every circle action of the form $s_i \pm n t_i$) is such that the number of full one-turns of that action at the fiber over $[0,1]$ minus the number of one-turns at the fiber over $[1,0]$ is equal to $i$. Thus the circle action $s_{2k-1}$ is mapped by $\phi$ to $s_{2k} - t_{2k}$. A similar argument shows that $s_{2k-1}$ is mapped to  $s_{2k-2} + t_{2k-2}$ in the second case. \end{proof}
 
\MS
We now prove the equivalence between configuration spaces and strata of almost complex structures.

\begin{prop}\label{Prop:EquivalenceJi&Ci}
The space $\tjj_i$ is homotopy equivalent to the configuration space ${\mathcal C}_i$.
\end{prop}
\begin{proof}
Consider the map $\pi:\tjj_i\to {\mathcal C}_i$ which sends each almost complex structure $\tJ$ to the triple of $\tJ$-curves in classes $F-E$, $E$ and $D_i$.  We first show
that the inverse image by $\pi$ of any configuration $C=(A_1,A_2,A_3)$ is non-empty and contractible. But $\pi^{-1}(C)$ consists of those $\tilde{\om}_{\mu,c}$-compatible structures $\tJ$ that  have all the tangent spaces to the configuration as invariant subspaces, and which are positive on the oriented surfaces $A_1$, $A_2$, and $A_3$. This space fibers over the contractible space $\Jj(C, \tilde{\om}_{\mu,c} |_C)$ of positive almost complex structures defined only on the union of the three curves.  The fiber of $\pi^{-1}(C)\to\Jj(C, \tilde{\om}_{\mu,c} |_C)$ retracts to the subspace of structures that are fixed on the normal directions to $C$ and therefore to the subspace of those almost complex structures which are fixed on some open neighborhood of $C$. But it is well known that the set of $\om$-tamed almost complex structures that are fixed on some open subset of a symplectic manifold is contractible.

To complete the proof, it is enough to show that $\pi$ is a fibration. Actually, since $\pi$ is a smooth map between Fr\'echet manifolds, it is sufficient to prove that $\pi$ is a submersion (see~\cite{Me}).

Choose some configuration $\Gamma$, and consider a tangent vector $v\in T_\Gamma\Cc_i$ represented by some path $\Gamma_t$. Given any tamed almost complex structure $\tJ\in\pi^{-1}(\Gamma)$, we can construct a vector $w\in T_{\tJ}\tjj_i$ such that $d\pi_{\tJ}(w)=v$ in the following way. The germ of the path $\Gamma_t$ is generated by a smooth isotopy $\phi_t\subset\Diff$. Since the tameness condition is open, there exists $\epsilon>0$ such that the path $J_t=\phi_t(J)$, $t<\epsilon$, lies in $\tjj_i$. By construction, this path $\tJ_t$ defines a vector $w$ such that $d\pi_{\tJ}(w)=v$.
\end{proof}
\BS
\begin{prop}\label{Prop:EquivalenceCi&Ci0}
The configuration space ${\mathcal C}_i$ is homotopy equivalent to the space ${\mathcal C}_i^0$ of orthogonal configurations.
\end{prop}
\begin{proof}
Since the three curves of each configuration intersect at exactly $2$ distinct points, we have two fibrations
$$
\begin{array}{ccccc}
{\mathcal C}_{p,q} & \to & {\mathcal C}_i & _\searrow & \\
\uparrow & & \uparrow & & \tilde{M}_{\mu,c} \times \tilde{M}_{\mu,c} - \mbox{~diag}\\
{\mathcal C}_{p,q}^0 & \to & {\mathcal C}_i^0 & ^\nearrow &
\end{array}
$$
where ${\mathcal C}_{p,q}$ is the space of configurations whose curves intersect at $p$ and $q$ and where the vertical maps are inclusions. Since these two fibrations have the same base,  one only needs to prove that the inclusion ${\mathcal C}_{p,q}^0\to {\mathcal C}_{p,q}$ is a homotopy equivalence. Considering now the tangent planes to the central curve $A_2$ at the points $p$  and $q$, we define two other fibrations
$$
\begin{array}{ccccc}
{\mathcal C}_{p,q,P,Q} & \to & {\mathcal C}_{p,q} & _\searrow & \\
\uparrow & & \uparrow & & G_{2,4}(T_p\tilde{M}_{\mu,c})\times G_{2,4}(T_q\tilde{M}_{\mu,c}) \\
{\mathcal C}_{p,q,P,Q}^0 & \to & {\mathcal C}_{p,q}^0 & ^\nearrow & \\
\end{array}
$$
where the base is a product of symplectic Grassmannian manifolds and where the fibers consist of configurations whose $A_2$-curve has fixed tangencies at $p$ and $q$. To show that the inclusion ${\mathcal C}_{p,q,P,Q} \to {\mathcal C}_{p,q,P,Q}^0$ is a weak homotopy equivalence, we can now work locally and choose two Darboux charts centered at $p$ and $q$ for which the images of $A_2$ are tangent to the vertical plane at the origin in $\R^2\times \R^2$.
We can then proceed as in~\cite{MP} Proposition 4.1.C to define local symplectic isotopies that map the images of $A_1$ and $A_3$ to curves orthogonal to $A_2$. Moreover, all these isotopies can be made to depend continuously on any compact family of parameters. This proves that the spaces ${\mathcal C}_{p,q,P,Q}$ and  ${\mathcal C}_{p,q,P,Q}^0$ are weakly homotopy equivalent.
\end{proof}
\BS
Consider the subgroup $\Symp_{ht}(\tMuc)\subset\Symp(\tMuc)$ of those symplectomorphisms that act trivially in homology. (Note that it coincides with the group $\Symp(\tMuc)$ when $\mu > 1$.) The next proposition shows that $\Cc^0_i$ is an homogeneous space under the action of $\Symp_{ht}(\tMuc)$.
\begin{prop}\label{Prop:Transitivite}
Assume that $\mu$ is rational. The natural action of $\Symp_{ht}(\tilde{M}_{\mu,c})$ on each configuration space ${\mathcal C}_i^0$ is transitive.
\end{prop}
\begin{proof}
Let $\om_0$ denote the standard symplectic form on $\tilde{M}_{\mu,c}$, and consider two orthogonal configurations $C_0$ and $C_1$ of the same type. By the symplectic neighborhood theorem, we can find a neighborhood ${\mathcal U}$ of $C_0$ and a diffeomorphism $\phi:(\tilde{M}_{\mu,c},C_0)\to (\tilde{M}_{\mu,c},C_1)$ whose restriction to ${\mathcal U}$ is a symplectomorphism. Consider the pull-back form $\om_1=\phi^*(\om_0)$. By construction, this form agrees with $\om_0$ near $C_0$; hence its restriction to $\tilde{M}_{\mu,c}-C_0$ is standard near the boundary. By Lemma~\ref{le:Liouville} which we postpone to Section~4.6, there is a Liouville flow $\psi_t$ that retracts $\tilde{M}_{\mu,c}-C_0$ into some standard polydisc $P$. After renormalization, this flow carries $\om_1$ to a form $\Om_1$ which is standard near the boundary of $P$. We can then apply Gromov's theorem and find a diffeomorphism $\xi$, isotopic to the identity, whose support is contained in $P$ and such that $\xi^*(\Om_1)$ is the standard form. Thus conjugating $\xi$ with the flow $\psi_t$ defines a symplectomorphism $(\tilde{M}_{\mu,c}-C_0,\om_0)\to (\tilde{M}_{\mu,c}-C_0,\om_1)$ that is the identity near the boundary. Composing this diffeomorphism with $\phi$, we get a symplectomorphism of $(\tilde{M}_{\mu,c},\om_0)$ that maps $C_0$ onto $C_1$. Since the homology classes represented by the curves of $C_i$ form a basis of $H_2(\tilde{M}_{\mu,c})$, this symplectomorphism induces the identity in homology.
\end{proof}
\BS
\begin{prop}\label{Prop:EquivalenceHi&Ti}
For $\mu$ rational, each stabilizer $\tH_i$ of the action of $\Symp_{ht}(\tilde{M}_{\mu,c})$ on ${\mathcal C}_i^0$  is homotopy equivalent to the corresponding torus $\tT^2_i$ of  K\"ahler isometries.
\end{prop}
\begin{proof}
The evaluation map
$$\tH_i \to \Symp(A_1,p)\times\Symp(A_2,\{p,q\})\times\Symp(A_3,q)$$
defines a locally trivial fibration whose fiber is the group $\tH_\id$ of symplectomorphisms  that fix the configuration $C_i$ pointwise. Since the base  retracts onto its subgroup $S^1\times S^1\times S^1$, the stabilizer $\tH_i$ is homotopically  equivalent to the subgroup $\tH_i'$ of symplectomorphisms acting linearly on $C_i$. Note that, by Proposition~\ref{prop:action}, the action of $\tT^2_i$ restricted to the curves $A_1$ and $A_2$ generates $S^1\times S^1\subset\Symp(A_1,p)\times\Symp(A_2,p)$. Thus the quotient $\tH'_i/\tT^2_i$ is homeomorphic to the subgroup
$$
\tH_i''=\{\phi\in\tH_i' ~~|~~ \phi|_{A_1\cup A_2}=\id,~~\phi|_{A_3}\mbox{~is linear}\}
$$
By the symplectic neighborhood theorem, the behavior of a symplectomorphism $\phi\in\tH_i''$  in a neighborhood of $A_1\cup A_2$ is determined, up to isotopy, by the action of its differential $d\phi$ on the symplectic normal bundles $\nu_1=\nu(A_1)$, $\nu_2=\nu(A_2)$. The restriction of the differential to $\nu_1$ and $\nu_2$ defines a fibration
$$
\begin{array}{ccccc}
\tH_i''' & \to & \tH_i'' & \to     & \Aut_p(\nu_1)\times\Aut_p(\nu_2) \\
        &     & \phi  & \mapsto & (d\phi|_{A_1},d\phi|_{A_2})
\end{array}
$$
with typical fiber $\tH_i''' = \{\phi\in\tH_i'' ~|~ \phi=\id\mbox{~on $\nu_1$ and $\nu_2$}\}$. Up to homotopy, $\nu_1$ and $\nu_2$ are just principal $U(1)$-bundles.  But the groups $\Aut_p(\nu_k)\simeq C^\infty((S^2,*);(S^1,1))$, $k=1,2$ because $d\phi_p=\id$, and are therefore contractible. The space $\tH_i''$ is thus homotopy equivalent to $\tH_i'''$. Observe that $d\phi_q=\id$ for all $\phi\in\tH_i'''$. Since each $\phi\in\tH_i'''$ acts linearly on $A_3$, its restriction to the curve $A_3$ is the identity. Hence, there is a fibration
$$
\begin{array}{ccccc}
 \tH_i''' & \to     & \Aut_q(\nu_3)  & \simeq & \{*\}\\
 \phi     & \mapsto & d\phi|_{\nu_3} &        &
\end{array}
$$
Its fiber is homotopy equivalent to the group of symplectomorphisms equal to the identity near $C_i$. But this group is contractible by Lemma~\ref{le:Liouville} (see the end of this section). Therefore, the groups $\tH_i'''\simeq\tH_i''\simeq \tH'_i/\tT^2_i$ are all contractible.
\end{proof}
\begin{cor}\label{Cor:HomogeneiteDesStrates}
For $\mu$ rational, the space $\tjj_i$ is homotopy equivalent to the homogeneous space $\Symp_{ht}(\tilde{M}_{\mu,c}) / \tT^2_i$.
\end{cor}
\begin{cor}\label{Cor:EquivalenceG&T2}
When $\mu = 1$, the identity component $\Symp_0(\tMuc)\subset\Symp(\tMuc)$ is homotopy equivalent to a torus $T^2$, while the full group is homotopy equivalent to the obvious extension of $T^2$ by $\Z_2$.
\end{cor}
\begin{proof}
When $\mu = 1$, there is only one stratum, which is therefore contractible, namely,
$$\tjj(\omlc)=\tjj_0\simeq \Symp_{ht}(\tilde{M}_{\mu,c}) / \tT^2_0$$
which implies the homotopy equivalence of $\Symp_{ht}(\tilde{M}_{\mu,c})$ with its connected subgroup $\tT^2_0$ and hence with $\Symp_0(\tMuc)$. The second assertion follows from the fact that the two $S^2$ factors of $S^2\times S^2$ can be interchanged by a symplectomorphism fixing the standard embedded ball $B_c$ and from the obvious fact that no other nontrivial automorphism of $H_2(\tilde{M}_{\mu,c})$ can be realized by a symplectic diffeomorphism.
\end{proof}
\begin{prop}\label{Prop:InclusionEnHomologie}
The inclusion $\tilde{T^2_i}\hookrightarrow \Symp(\tMuc)$ induces an injection in rational homology for all $i \ge 0$ and $\mu \ge 1$.
\end{prop}
\begin{proof}
First note that the rational homology rings of $\tilde{T^2_i}$ and $\Symp(\tMuc)$ are both Hopf algebras. Since the inclusion $\tT^2_i\hookrightarrow \Symp(\tMuc)$ induces a morphism of algebras, it is enough to show that the submodule of primitive elements of $H_*(\tT^2_i;\Q)$, which is simply $H_1(\tT^2_i;\Q) 
$, is mapped injectively into $H_1(\Symp(\tMuc);\Q)=\pi_1(\Symp(\tMuc))\otimes\Q$.

Consider first the case $i=2k$. The blow-down of the exceptional sphere in class $E$ induces a continuous map $\Symp(\tMuc,E) \to \Ee(W_{2k},*)$, where $\Ee(W_{2k},*)$ is the space of  homotopy self-equivalences of the Hirzebruch surface $W_{2k}\simeq S^2\times S^2$, preserving the point~$*$ equal to the intersection of the zero section with the fiber. This map sends $\tilde{T^2_{i}}$ onto the torus $T^2_{2k}$. Therefore it is sufficient to show that the natural map $\pi_1(T^2_{2k})\otimes\Q\to\pi_1(\Ee(W_{2k},*))\otimes\Q$ is injective.

The action of $T_{2k}$ on $W_{2k }$ is given by
\begin{eqnarray}\label{Action}
(\theta_1,\theta_2)\cdot([z_0,z_1],[w_0,w_1,w_2]) & = & ([\theta_1z_0,z_1],[\theta_1^{2k}w_0,w_1,\theta_1^{2k}\theta_2w_2])
\end{eqnarray}
We can use the zero section $s_0:[z_0,z_1]\mapsto([z_0,z_1],[0,0,1])$ and the projection to the first factor $\pi_1:([z_0,z_1],[w_0,w_1,w_2])\mapsto [z_0,z_1]$ to define a continuous map
\begin{eqnarray*}
\Ee(W_{2k},*) & \to     & \Ee(S^2,*) \\
       \phi & \mapsto & \pi_1 \circ \phi \circ s_0
\end{eqnarray*}
Let $s_{2k}$ and $t_{2k}$ be the generators of $\pi_1(T^2_{2k})$. From equation~(\ref{Action}) it is clear that the induced map $\pi_1(\Ee(W_{2k},*))\to\pi_1(\Ee(S^2,*))$ sends $s_{2k}$ to an element of infinite order  while $t_{2k}$ is sent to $0$. Using the inclusion of $S^2$ as the fiber $f:[z_0,z_1]\mapsto([1,0],[z_0,0,z_1])$ and the projection to the second factor $\pi_2:W_{2k}\to S^2$, one may define a  map
\begin{eqnarray*}
\Ee(W_{2k},*) & \to     & \Ee(S^2,*) \\  
       \phi & \mapsto & \pi_2 \circ \phi \circ f
\end{eqnarray*}
which sends $t_{2k}$ to an element of infinite order in $\pi_1(\Ee(S^2,*))$.
Therefore, both $s_{2k}$ and $t_{2k}$ are of infinite order in $\pi_1(\Ee(W_{2k},*))$, and there is a map sending $t_{2k}$ (but not $s_{2k}$) to zero. Thus any rational relation $ms_{2k}+nt_{2k}=0$ in $\pi_1(\Ee(W_{2k},*))$ must be trivial.

Consider now the odd case $i=2k-1$. By the relations in Proposition~\ref{prop:invert}, the actions of $\tT^2_{2k-1}= (\tT^2_{2k-1})_0$ and of  $(\tT^2_{2k})_{\infty}$ on $\tMuc$ coincide. Thus the blow-down of this action along the class $E=[A_1] \in H_2(\tMuc)$ is equal to the 2-torus of K\"ahler isometries of $W_{2k}$ that fix the point $s_{\infty} \cap F_i$, and we may conclude using the same argument as in the above even case. 

   There is another, perhaps simpler, way of proving this proposition. Indeed, the same proof as in Lemma~\ref{le:retracts} shows that the space $\Symp(\tMuc)$ is homotopy equivalent to $\Symp(\tMuc, E \cup (F-E))$ since none of these curves can degenerate. One therefore gets a map from $\Symp(\tMuc) \to \Symp(\tMuc, E \cup (F-E)) \to \Ee(W_{i}, F_i, *)$, the space of homotopy self-equivalences of $W_i$ that preserve $F_i$, and $*$, the point of intersection between $F_i$ and the section of self-intersection $-i$. Here the last map is the one induced by blowing down the $A_2$-curve which is $E$ in the even case and  $F-E$ in the odd case. Consider the maps
$$
\Ee(W_{i}, F_i, *) \to \Ee(F_i, *)
$$
given by the restriction to the fiber $F_i$ and 
$$
\Ee(W_{i}, F_i, *) \to \Ee(S^2, *)
$$
given as above by $ \phi  \mapsto  \pi_1 \circ \phi \circ s_i$ (where here $S^2$ is the base of the Hirzebruch fibration). They obviously show that the two elements $s_i,t_i$ are independent.
\end{proof}
\begin{prop}\label{Prop:ScindementHomologique}
For any rational $\mu\geq 1$ and any $0<c<1$,
\begin{enumerate}
\item[(1)] the group $\Symp_{ht}(\tMuc)$ is connected;
\item[(2)] for all $0\leq i \leq \ell$, the stratum $\tjj_i$ is connected;
\item[(3)] for all $0\leq i \leq \ell$, $H_*(\Symp_{ht}(\tilde{M}_{\mu,c})) \simeq H_*(\tjj_i)\otimes H_*(\tT^2_i)$ as graded vector  spaces and $H^*(\Symp_{ht}(\tilde{M}_{\mu,c})) \simeq H^*(\tjj_i)\otimes H^*(\tT^2_i)$ as graded algebras.
\end{enumerate}
\end{prop}
\begin{proof}
Consider the fibration
$$
\Symp_{ht}(\tilde{M}_{\mu,c})\to\Symp_{ht}(\tilde{M}_{\mu,c})/\tT^2_i\simeq \tjj_i
$$
Its fiber is obviously connected and finite-dimensional. When $i=0$, the base $\tjj_0\simeq\Symp_{ht}(\tMuc)/\tT^2_0$ is also connected, and we conclude at once that $\Symp_{ht}(\tMuc)$ is itself connected. In turn, the connexity of $\Symp_{ht}(\tMuc)$ implies that each stratum $\tjj_i$ is connected for all $i\leq\ell$. By Proposition~\ref{Prop:ScindementHomologique}, the fiber is totally nonhomologous to zero. Therefore, we can apply the Leray-Hirsch theorem to conclude that
$$H_*(\Symp_{ht}(\tilde{M}_{\mu,c})) \simeq H_*(\tjj_i)\otimes H_*(\tT^2_i)$$
as vector spaces.

The same argument shows that, as graded vector spaces,
$H^*(\Symp_{ht}(\tilde{M}_{\mu,c})) \simeq H^*(\tjj_i)\otimes H^*(\tT^2_i)$. Moreover, since $H^*(\tT^2_i)$ is a free algebra,  we can define a homomorphism of algebras $\theta:H^*(\tT^2_i) \to H^*(\Symp_{ht}(\tilde{M}_{\mu,c}))$ which is a right inverse to the restriction. This shows that the above identification is an isomorphism of algebras.
\end{proof}
\begin{cor}\label{Cor:EquivalenceHomologiqueDesStrates}
Assuming $\mu$ rational, $H_*(\tjj_i)\simeq H_*(\tjj_0)$ for all $1\leq i \leq \ell$.
\end{cor}
\begin{remark}
We note here that although $\Symp(\tilde{M}_{\mu,c})$ is a topological group, one cannot assume a priori that its rational cohomology is a Hopf algebra. A sufficient condition for this to be true is that its rational homology $H_*(\Symp(\tilde{M}_{\mu,c}))$ is of finite type, that is, that each vector space $H_q(\Symp(\tilde{M}_{\mu,c}))$ is finite-dimensional. By the previous proposition, this is equivalent to $H_*(\tjj_i)$ being of finite type for one (and hence any) $0\leq i \leq \ell$. We show in Sections~\ref{BigCase} and~\ref{SmallCase} that these homologies are indeed of finite type.
\end{remark}

What we did in this section shows that, as the parameter $c$ crosses the critical value $\mu-[\mu]_-$, the number of strata jumps in the space $\tjj(\omlc)$ of almost complex structures on $\tMuc$ but this has of course no effect on the stratification of $\jj(\om_{\mu})$ on the manifold $M_{\mu}$. Thus we may expect that the cohomology of $\Symp(\tilde{M}_{\mu,c})$ jumps at that value, whereas the cohomology of $G_{\mu}$ remains clearly
unchanged. This suggests that the space
$\Im \Emb_{\om}(c,\mu)= \Im \Emb_{\om}(B^4(c),M_{\mu})$, which is the homotopy quotient of these two groups, undergoes a homotopy jump at that value. This is indeed the case in general (see \cite{Pinso}). In the next sections, we will establish this for the case  $\mu \le 2$.

\subsection{The case $1<\mu\leq 2$ and $0 < \lambda \leq c < 1$}\label{BigCase}

When $0 < \lambda \leq c < 1$, the space $\tjj$ is the disjoint union of $2$ strata.
\begin{lemma}\label{Lemme:Isomorphisme2Strates}
When $0 < \lambda \leq c < 1$, there is an isomorphism
$H_{i+1}(\tjj_0)\simeq H_i(\tjj_1)$, $\forall\, i\geq 0$.
\end{lemma}
\begin{proof}
Since the space $\tjj_{01} = \tjj(\om)$ is contractible, the long exact sequence of the pair  $(\tjj_{01},\tjj_0)$ shows that the groups $H_{i+2}(\tjj_{01},\tjj_0)$ and $H_{i+1}(\tjj_0)$ are isomorphic. Let $N(\tjj_1)$ be a neighborhood of $\tjj_1$ in the contractible space $\tjj_{01}$, and write $\xi$ for the intersection $\tjj_0\cap N(\tjj_1))$. Excising $\tjj_{01}-N(\tjj_1)$ from $(\tjj_{01},\tjj_0)$, we see that   $H_{i+2}(\tjj_{01},\tjj_0)\simeq H_{i+2}(N(\tjj_1), \xi)$. Since $\tjj_1$ is a codimension $2$ co-oriented and connected submanifold of $\tjj_{01}$, its neighborhood $N(\tjj_1)$ is homeomorphic to an oriented fiber bundle of rank 2 over $\tjj_1$. Hence, by the Thom isomorphism theorem, there is an isomorphism
$H_{i+2}(N(\tjj_1), \xi) \simeq H_{i}(\tjj_1)$ for all $i\geq 0$.
\end{proof}
\begin{cor}\label{Cor:Homologie2Strates}
When $0 < \lambda \leq c < 1$, the rational homology groups of $\Symp(\tilde{M}_{\mu,c})$ are
$$
H_q(\Symp(\tilde{M}_{\mu,c})) =
\left\{
\begin{array}{ll}
   \qq    & \mbox{if $q = 0$}     \\
   \qq^3  & \mbox{if $q = 1$}     \\
   \qq^4  & \mbox{if $q \geq 2 $} \\
\end{array}
\right.
$$
In particular, the graded module $H_*(\Symp(\tilde{M}_{\mu,c}))$ is of finite type.
\end{cor}
\begin{proof}
By Theorem~\ref{Thm:Stabilite}, we may assume that $\mu$ is rational and can therefore apply all statements proved so far. Lemma~\ref{Lemme:Isomorphisme2Strates} and Corollary~\ref{Cor:EquivalenceHomologiqueDesStrates} together imply that $H_{q+1}(\tjj_1)\simeq H_{q+1}(\tjj_0)\simeq H_q(\tjj_1)$, for all $q\geq 0$. Since each stratum is connected, $H_0(\tjj_1) \simeq \Q$. Hence $H_i(\tjj_1) \simeq \Q$ for all $i \ge 0$. The result now follows from the fact that
$H_*(\Symp(\tMuc))\simeq H_*(\tjj_1)\otimes H_*(\tilde{T}^2_1)$.
\end{proof}
Being of finite type, the rational cohomology of $\Symp(\tilde{M}_{\mu,c})$ is an associative and commutative Hopf  algebra dual to $H_*(\Symp(\tilde{M}_{\mu,c}))$. By Leray's structure theorem, it is isomorphic, as an algebra, to the product of an exterior algebra generated by elements of odd degrees with a symmetric algebra generated by elements of even degrees. It follows from the previous corollary that
\begin{prop}\label{StructureCohomologie2Strates}
When $0 < \lambda \leq c < 1$, the rational cohomology of $\Symp(\tilde{M}_{\mu,c})$ is isomorphic, as an algebra,
to the product
$$\Lambda(\alpha_1,\alpha_2,\alpha_3)\otimes S(\epsilon)$$
where each $\alpha_i$ is of degree $1$ and where $\epsilon$ is of degree $2$.
\end{prop}
We compute in Section~\ref{se:Pontrjagin} the rational homology ring, that gives more informations on the algebraic structure of the group.

\subsection{The case $1<\mu\leq 2$ and $0<c<\lambda\leq 1$}\label{SmallCase}

When $0 < c < \lambda \leq 1$, the space $\tjj$ is the disjoint union of $3$ strata and we could proceed as above to calculate the rational cohomology ring of $\Symp(\tMuc)$. Instead, we use Theorem~\ref{Thm:Stabilite} to prove a stronger result, namely, the following.
\begin{prop}\label{StructureCohomologie3Strates}
When $0 < c < \lambda \leq 1$, the group $\Symp(\tMuc)$ has the same homotopy  type as the stabilizer of a point $\Symp_p(M_\mu) \subset \Symp(M_\mu)$. In particular, its rational cohomology is isomorphic, as an algebra, to the product
$$\Lambda(\alpha_1,\alpha_2,\alpha_3)\otimes S(\epsilon)$$
where each $\alpha_i$ is of degree $1$ and where $\epsilon$ is of degree~$4$.
\end{prop}
\begin{proof}
We first calculate the rational cohomology of $\Symp_p(M_\mu)$. By the work of Abreu~\cite{Ab}, the rational homotopy of $\Symp(M_\lambda)$ is
$$
\pi_q(\Symp(M_\lambda))\otimes \Q = \left\{
\begin{array}{ll}
\Q   & \mbox{~if $q=1$ or $q=4$}\\
\Q^2 & \mbox{~if $q=3$}\\
 0   & \mbox{~otherwise}
\end{array}
\right.
$$
where the two generators of degree $3$ correspond, respectively, to the standard action of $SO(3)$ on the first factor of $M_\mu = S^2\times S^2$ and the action of $SO(3)$ on the diagonal, while the generator of degree $1$ is the rotation in the fibers of $W_2$. Looking at the homotopy sequence of the fibration $\Symp_p(M_\mu) \to \Symp(M_\mu) \to M_\mu$, we see that the rational homotopy of the  stabilizer $\Symp_p(M_\mu)$ is
$$\pi_q(\Symp_p(M_\mu))\otimes \Q = \left\{
\begin{array}{ll}
\Q^3 & \mbox{~if $q=1$}\\
\Q   & \mbox{~if $q=4$}\\
 0   & \mbox{~otherwise}
\end{array}
\right.
$$
In particular, the rational homotopy of $\Symp_p(M_\mu)$ is finitely generated, and both graded modules $H_*(\Symp_p(M_\mu))$ and $H^*(\Symp_p(M_\mu))$ are of finite type. Therefore the classical Cartan-Serre theorem (see Section~\ref{se:Pontrjagin}) implies that the rational cohomology of $\Symp_p(M_\mu)$ is generated, as an algebra, by the duals of its spherical classes. This proves that $H^*(\Symp_p(M_\mu))$ is isomorphic to the algebra $\Lambda(\alpha_1,\alpha_2,\alpha_3)\otimes S(\epsilon)$, where each $\alpha_i$ is of degree $1$ and where $\epsilon$ is of degree~$4$.

We now prove that the stabilizer $\Symp_p(M_\mu)$ is homotopy equivalent to the group $\Symp^{\U(2)}(\tMuc,E)$. First, recall that every symplectomorphism of $\Symp^{\U(2)}(\tMuc,E)$ gives rise to a symplectomorphism of $M_\mu$ acting linearly near the embedded ball $B_{c}$ and fixing its center $p$. Conversely, every homotopy class of the stabilizer $\Symp_p(M_\mu)$ can be realized by a family of symplectomorphisms that act linearly on a ball $B_{c'}$ of sufficiently small capacity $c'$ centered at $p$. Hence, we can lift such a representative to the group $\Symp^{\U(2)}(\tM_{\mu,c'},E)$ of symplectomorphism acting linearly near $E$. Obviously the directed system of homotopy equivalences
$$
r_{c',c}: \Symp^{\U(2)}(\tM_{\mu,c'},E) \to \Symp^{\U(2)}(\tMuc,E)
$$
where $r$ is the restriction map, commutes with the maps
$$
g_{c}: \Symp^{\U(2)}(\tMuc,E) \to \Symp_p(M_\mu)
$$
and yields a commutative triangle for each pair $c \le c'$. It is then obvious that each map $g_c$ is a weak homotopy equivalence. Finally, compose this homotopy equivalence with the equivalence
$$
\Symp^{\U(2)}(\tMuc,E)  \to \Symp(\tMuc)
$$
\end{proof}

Denote by $a_1$ and $a_2$ the blow-up of the rotations in the first factor and the diagonal of $M_\mu$ and write $a_3$ for the blow-up of the $S^1$-action on $M_\mu$ that generates $\pi_1(\Symp(M_\mu))\otimes\Q$. Observe that the elements $a_1$ and $a_2$ generate $\pi_1(\tT_0^2)$ while $a_3$ can be identified with one generator of $\tT_1^2$. By the proof of Proposition~\ref{StructureCohomologie3Strates}, $e$ can be identified with the blow-up of the generator of $\pi_4(\Symp(M_\mu))\otimes\Q$. The last proposition implies the following.
\begin{cor}\label{GenerateursStabilisateur}
When $0<c < \lambda \le 1$, the rational homotopy of $\Symp(\tMuc)$ is generated, as a module, by the elements $a_1$, $a_2$, $a_3$ and $e$.
\end{cor}
%

\subsection{The Lie algebra $\pi_*(\Symp(\tMuc))\otimes\Q$ and the
rational Pontrjagin ring of $\Symp(\tMuc)$}\label{se:Pontrjagin}

The Milnor-Moore theorem on Hopf algebras states that a connected cocommutative (in the graded sense) Hopf algebra  $A$ on a field of characteristic zero is generated by its primitive elements. A {\em primitive element} is an element  $a \in A$ such that its coproduct is equal to $1 \otimes a + a \otimes 1$. More precisely, the  Milnor-Moore theorem says that any such Hopf algebra is isomorphic, as Hopf algebra, to the enveloping algebra ${\mathcal E}$ of the Lie subalgebra consisting of the primitive elements. This means more precisely that the only relations in ${\mathcal E}$  are the ones given by $a\otimes b - (-1)^{pq}b \otimes a = a b - (-1)^{pq}b a$ where the product on the right is the product in the Hopf algebra.

When $A$ is the rational homology of an H-space, the coproduct is induced by the inclusion of  $H(G)$ in the diagonal $H(G\times G) \simeq H(G) \otimes H(G)$, so  a primitive element is an element such that  its inclusion in $H(G \times G)$, via the diagonal map, decomposes by the Kunneth  formula as $1\otimes a+ a\otimes 1$. The Cartan-Serre theorem states that the primitive elements, in the case of the homology of a H-space with coefficients in a field of characteristic zero, are precisely the spherical classes. The Lie product of two classes $a$, $b$ is then the Samelson product $[a,b]$ which, in the case of a topological group $G$, is induced by the commutator
\begin{eqnarray*}
S^p\times S^q & \to     & G \\
     (u,v)~~  & \mapsto & a(u)b(v)a^{-1}(u)b^{-1}(v)
\end{eqnarray*}
Therefore we have the following (see~\cite{CS}).
\begin{thm} [Cartan-Serre, Milnor-Moore] \label{thm:MM} Let $K$ be a $H$-space, and consider its rational homology, with the Pontrjagin product induced by the product in $K$ and the coproduct induced by the inclusion of the diagonal. Let $\Uu(\pi_*(K)\otimes Q)$ be the enveloping algebra with the only relations given by $a\otimes b - (-1)^{pq}b \otimes a = [a,b]$ where the bracket is the Samelson product. Then the natural map induced by the Hurewicz homomorphism
$$
\Uu(\pi_*(K)\otimes Q) \to H_*(K; \Q)
$$
is an isomorphism of Hopf algebras.

Moreover, if the rational homology is finitely generated in each dimension, the rational cohomology is a Hopf algebra for the cup product and coproduct induced by the product in $H$, and it is generated as an algebra by elements that are dual to the spherical classes in homology. In particular, the number of generators of odd dimension $d$ appearing in the antisymmetric part $\Lambda$ is equal to the dimension of $\pi_d(G) \otimes \Q$, and the number of generators of even dimension $d$ appearing in the symmetric part $S$ is equal to the dimension of $\pi_d(G) \otimes \Q$.
\end{thm}

\MS
We now compute the Lie algebra $\pi_*(\Symp(\tMuc))\otimes\Q$ for $\mu\le 2$. For this, we need to describe the generators of this algebra in more detail when $c \ge \la$. In this case, there are two strata and therefore two tori $\tilde T^2_0$ and $\tilde T^2_1$. The first one is generated by the elements $x_0, y_0 \in \pi_1(\tilde T^2_0)$, where $x_0, y_0$ are the blow-ups of the rotations in the first  and second factors of $S^2\times S^2$ that fix the coordinate $z$ of the point $p=(z,z)$. Set $a_1 = x_0$ and $a_2 = x_0 + y_0$ the diagonal rotation. The second torus $\tilde T^2_1$ is generated by the elements $x_1,y_1$ which are blow-ups of the $S^1$-rotation $s_1,t_1$ described in Proposition~\ref{prop:action}. By Proposition~\ref{prop:invert}, the generator $x_1$ of the torus $(\tilde{T}^2_1)_0$ is equal to the generator $(x_2)_\infty - (y_2)_\infty$ of $(\tilde{T}^2_2)_\infty$. Hence, $x_1+y_1$ corresponds to $(x_2)_\infty$. But $(x_2)_\infty$ is, by definition, the blow-up at point $p$ of the generator $s_2$, that is, the generator of the $\pi_1$ of the $\SO(3)$-part of the K\"ahler isometries of $W_2\simeq S^2\times S^2$ which fixes $p$. We know by Abreu's work that this $\SO(3)$-action may be identified to the standard diagonal action of $\SO(3)$ on the product $S^2\times S^2$. Hence, $(x_2)_\infty$ corresponds to the blow-up of the $S^1$-diagonal in the torus $T_0=S^1\times S^1\subset \SO(3)\times \SO(3)$. That is to say, $(x_2)_\infty$ is identified with $x_0+y_0$. Therefore, in the basis $\{x_1+y_1, y_1\}$ of $\tilde{T}^2_1$, the first element can be identified with $a_2=x_0+y_0$. It is clear that the elements $a_1=x_0,a_2 = x_0 + y_0,y_1$ generate $H_1(\Symp(\tMuc))$. Moreover, we know that all these generators commute except perhaps the elements $y_1$ and $x_0$. Actually, we have the following.
\begin{lemma} \label{NotCommuting}
The elements $y_1$ and $x_0$ do not commute; that is, their Samelson product does not vanish.
\end{lemma}
\begin{proof}
Consider the evaluation maps
\begin{eqnarray*}
p^i:\Symp(\tMuc) & \to & \tjj_i\simeq\Symp(\tMuc)/\tT_i \\
\phi & \mapsto & \phi\cdot\tJ_i
\end{eqnarray*}
It is sufficient to show that the Samelson product $[x_0,y_1]$ maps to a nonzero element in $H_2(\tjj_0)$. The Hurewicz homomorphism $h_2$ sends $[x_0,y_1]$ to $x_0y_1+y_1x_0$ in the Pontrjagin ring $H_*(\Symp(\tMuc))$. Since $x_0$ stabilizes the complex structure $\tJ_0$, the image of $y_1x_0$ in $\tjj_0$ is the one dimensional cycle $p^0(y_1)$. Therefore, $[x_0,y_1]$ is mapped to the cycle $p^0_*(x_0y_1)$ in $H_2(\tjj_0)$. We claim that this cycle is mapped to a nonzero element by the isomorphism
$$\psi:H_2(\tjj_0)\simeq H_1(\tjj_1)$$
from Lemma~\ref{Lemme:Isomorphisme2Strates}. Indeed, since $H_1(\tjj_0)$ is generated by the fundamental class $L$ of the link $\Ll = N(\tjj_1)\cap\tjj_0$, this isomorphism sends a class $p^0_*(z)$ to $p^1_*((p^0)^*(L^*)\cap z)$, where $L^*$ is dual to $L$. Now consider the isomorphism
$$H_1(\Symp(\tMuc)) \simeq H_1(\tjj_0)\oplus H_1(\tT^2_0)\simeq \Q\oplus\Q^2$$
Since $x_0$ and $y_0$ are the generators of $H_1(\tT^2_0)$, the element $p^0_*(y_1)\in H_1(\tjj_0)$ must represent $\pm L$. Therefore, $\psi(p^0_*(z))$ is equal to $p^1_*(y_1^*\cap z)$. Using the fact that both $x_0$ and $y_1$ are primitive elements of $H_*(\Symp(\tMuc))$, a direct calculation shows that  $y_1^*\cap (x_0y_1) = x_0$ and hence that $\psi(p^0_*(x_0y_1))=p^1_*(x_0)$.
Looking at the isomorphism
$$H_1(\Symp(\tMuc)) \simeq H_1(\tjj_1)\oplus H_1(\tT^2_1)\simeq \Q\oplus\Q^2$$
we see that $p^1_*(y_0)=p^1_*(x_1)=p^1_*(y_1)=0$ and hence that $p^1_*(x_0)$ must be non-zero. Therefore the class $p^0_*(x_0y_1)$ does not vanish.
\end{proof}
\begin{cor}\label{GenerateursHomotopie2Strates}
When $c \ge \la$, the elements $a_1=x_0, a_2=x_0 + y_0$ and $a_3=y_1$ generate the Pontrjagin ring of~$\Symp(\tMuc)$.
\end{cor}
\begin{proof}
This follows from Theorem~\ref{thm:MM}.
\end{proof}
\begin{thm}\label{thm:homotopy-tildeM}
\begin{enumerate}
\item[(1)] The Lie algebra $\pi_*(\Symp(\tMuc)) \otimes \Q$, with $\mu \in (1,2]$ and  $0 < c < \la$, is the commutative algebra generated by the elements $a_1,a_2,a_3, e$ of Corollary~\ref{GenerateursStabilisateur} and dual to the elements of Theorem~\ref{thm:cohomology}, with vanishing Samelson products.
\item[(2)] The Lie algebra $\pi_*(\Symp(\tMuc)) \otimes \Q$, with $\mu \in (1,2]$ and  $\la \le c < 1$, is generated by the elements $a_1,a_2,a_3, e$ dual to the elements of Theorem~\ref{thm:cohomology}. All Samelson products vanish except $[a_1, a_3] = \mbox{const~}e$.
\end{enumerate}
\end{thm}
\begin{proof}
(1) Since $\pi_q(\Symp(\tMuc)) \otimes \Q = 0$ except when $q=1$ or $q=4$, all Samelson products vanish.

(2) By the dual version of the Cartan-Serre theorem and Proposition~\ref{StructureCohomologie2Strates}, there is no rational homotopy sphere in dimension higher that $2$; hence the Samelson product of $e$ with any $a_i$ must vanish. By Lemma~\ref{NotCommuting}, the only nonvanishing Samelson product is $[a_1, a_3]$, which must therefore be equal to a nonvanishing multiple of $e$.
\end{proof}

For later use, we state the similar theorem for $\pi_*(\Symp(M_{\mu}))\otimes \Q$ proved by Abreu and McDuff in \cite{AM}:
\begin{thm} \label{thm:homotopy-M}
The Lie algebra $\pi_*(\Symp(M_{\mu})) \otimes \Q$ with $\mu \in (1,2]$ is generated by elements $a, b_1, b_2, d$ of degree $1,3,3$, and $4$, respectively, with vanishing Samelson products except for  $[a,b_1] = {\rm const} \, d$.
\end{thm}
The element $a$ is the rotation in the fibers of the Hirzebruch surface $W_2$, the element $b_1$ is the $3$-cycle of $SO(3)$ considered as acting on the first factor of $S^2 \times S^2$, and the element $b_2$ is the $3$-cycle of the diagonal action of $SO(3)$ on $S^2\times S^2$.

\MS
Theorem~\ref{thm:homotopy-tildeM}, together with Theorem~\ref{thm:MM}, gives at once the following.
\begin{thm}\label{thm:PontrjaginEclatement}
\begin{enumerate}
\item[(1)] The rational Pontrjagin ring of $\Symp(\tMuc)$, with $\mu\in (1,2]$ and $0<c<\la$, is isomorphic to $\Lambda(a_1,a_2,a_3)\otimes S(e)$; that is, it is isomorphic to its cohomology~ring.
\item[(2)] The rational Pontrjagin ring of $\Symp(\tMuc)$, with $\mu \in (1,2]$ and  $\la\le c<1$, is isomorphic to
$(\Q\langle a_1,a_2,a_3\rangle\otimes S(e) )/R$, where $R$ is generated by the relations $a_1^2=a_2^2=a_3^2=0$, $a_1a_2=-a_2a_1$, $a_2a_3=-a_3a_2$ and $a_1a_3+a_3a_1={\rm const~}e$.
\end{enumerate}
\end{thm}
%

\subsection{The inclusion $\Symp(\tMuc)\hookrightarrow\Diff(\tM)$}

When $c<\lambda$, the generator $e$ of the group $\pi_4(\Symp(\tMuc))$ corresponds, via the identification
$$\Symp(\tMuc)\simeq \Symp_p(M_\mu)\subset\Symp(M_\mu)$$
of Proposition~\ref{StructureCohomologie3Strates}, to the generator $d$ of $\pi_4(\Symp(M_\mu))$. As explained in~\cite{AM}, this element $d$ is in the kernel of the map $\pi_4(\Symp(M_\mu))\to\pi_4(\Diff(M))$. Hence looking at the long homotopy sequences of the fibrations
$$
\begin{array}{rrcccccll}
 \Symp(\tMuc)& \simeq&\Symp_p(M_\mu) & \to & \Symp(M_\mu) & \to & M_\mu & &\\
 & &\downarrow     &     & \downarrow   &     & \downarrow & &\\
 & &\Diff_p(M)     & \to & \Diff(M)     & \to & M &=& S^2\times S^2
\end{array}
$$
we see that $e$ maps to $0$ in $\pi_4(\Diff_p(M))$. Since we can lift the contracting homotopy in $\Diff_p(M)$ to $\Diff(\tM)$, we conclude that $e$ vanishes in $\pi_4(\Diff(\tM))$.

When $c\geq\lambda$, the generator $e$ of $\pi_2(\Symp(\tMuc))$ does not correspond to any element in $\pi_2(\Symp(M_\mu))$. On the other hand, if we identify the group $\Symp(\tMuc)$ with $\Symp^{\U(2)}(M_\mu,B_c)$, then for any $0< c_0 < c$, the restriction $\Symp^{\U(2)}(M_\mu,B_c)\to\Symp^{\U(2)}(M_\mu,B_{c_0})$ defines an embedding
$\Symp(\tMuc)\to \Symp(\tM_{\mu,c_0})$. Since the identifications depend continuously on the capacity, this embedding is isotopic to the inclusion $\Symp(\tMuc)\subset\Diff(\tM)$. Hence the generator $e\in\pi_2(\Symp(\tMuc))$ must vanish in $\pi_2(\Diff(\tM))$ since $\pi_2(\Symp (\tM_{\mu,c_0}))=0$ when $c_0<\lambda$. 

Now let us prove that the subspace spanned by the three generators of the group $\pi_1(\Symp(\tMuc)) \otimes \Q$, represented by Hamiltonian $S^1$-actions, injects in $\pi_1(\Diff(\tM)) \otimes \Q$.  We first show that this 3-dimensional subspace injects in the space $\pi_1(\Ee^0(\tM,C')) \otimes \Q$, the rational fundamental group  of the space of homotopy self-equivalences of $\tM$ that
\begin{enumerate}
\item[(1)] are homotopic to the identity through homotopy self-equivalences and
\item[(2)] preserve the standard configuration $C'=(F-E) \cup  E$ (not necessarily pointwise).
\end{enumerate}
Suppose that there is a relation
$$
\ell x_0 + m y_0 + n y_1 = 0 
$$
with $\ell, m, n \in \Z$, where $x_0, y_0, y_1$ are the $S^1$-actions considered in $\pi_1(\Ee^0(\tM,C')) \otimes \Q$.  As in the proof of  Proposition~\ref{Prop:InclusionEnHomologie}, the blow-down of $E$ yields a similar relation amongst the same generators considered in $\Ee^0(S^2 \times S^2, F,*)$. Recall that  $x_0,y_0,y_1$ are then respectively the one-turn rotations round the base $B$, the fiber $F$ and the rotation in the fibers round the curves in classes $B+F$ and $B-F$. The projection on the base
$$
\Ee^0(S^2 \times S^2, F, *) \to \Ee^0(B=S^2,*)
$$
sends $x_0$ to the rotation round $*$ and the other two to the constant loop. Hence $\ell = 0$. The restriction to the fiber $F$ shows that each of $y_0,y_1$ is of infinite order.  Consider the ``diagonal map''
$$
\Phi: \Ee^0(S^2 \times S^2, F, *) \to \Ee^0(S^2,*)
$$
defined by $\Phi(f) = \pi_2 \circ f \circ \si_{B+F}$. It sends each $f$ to the map from the base to the fiber defined by appying $f$ to the graph $B+F$ of the identity map $B \to F$. It is clear that $\Phi$ sends $y_1$ to the constant loop based at the identity map and $y_0$ to the generator of $\pi_1(\Ee^0(S^2,*)) \otimes \Q$. Thus $m=0$. But $n$ must also vanish because $y_1$ is of infinite order.

Consider now the fibration
$$
\Ee^0(\tM,C') \to \Ee^0(\tM) \to \Im \Map^0(C', \tM)
$$
where the second term is the space of homotopy self-equivalences which are homotopic to the identity through homotopy self-equivalences, and the last term is the space of images of smooth maps from $C'$ to $\tM$ which are homotopic to the standard inclusion. To show that $\Ee^0(\tM,C') \to \Ee^0(\tM)$ induces an injection at the rational $\pi_1$ level, it is enough to show that the map 
$$
\pi_2(\Im \Map^0 (C', \tM)) \otimes \Q \to  \pi_1(\Ee^0(\tM,C')) \otimes \Q
$$
is zero, and it is therefore sufficient to prove that the 
fibration $\Ee^0(\tM) \to \Im \Map^0(C', \tM)$ admits a lift at the rational $\pi_2$-level. Let  $f$ be a $S^2$-family of elements in $\Im \Map^0 (C', \tM)$ based at the image of the inclusion. Lift this family to a family $g$  of parametrised maps. By symmetry, there is no obstruction to extend this $S^2$-family to a 
$S^2$-family of maps
$$
\tilde g: S^2 \to Map^0 (C, \tM)
$$
based at the standard inclusion, where $C$ denotes the union of the standard representatives in classes $F-E,E,B-E$.  The complement of the standard configuration $C$ is the $4$-ball; denote by $h: S^3 \to C$ the attaching map.  For each element $z$ in the family, the composition $\tilde g_z \circ h$ is then contractible, and any extension to the 4-ball descends to a continuous map $\psi_z: \tM \to \tM$ which extends $\tilde g_z$ and is a homotopy self-equivalence. The extension from $\tilde g_z$ to $\psi_z$ can be made to depend continuously on $z$; actually the obstruction to such a continuous extension lies in the fundamental group of the space of self-homotopy equivalences of the 4-ball which restrict to the identity near the boundary. But this space is convex.

   This shows that the map $\Symp(\tMuc) \to \Ee^0(\tM)$ induces an injection at the rational $\pi_1$-level. But this map evidently factorizes through $\Diff^0(\tM)$. Therefore we have the following.

\begin{prop}
For any $\mu \in (1,2]$ and $c \in (0,1)$, the kernel of the inclusion\\
$
\pi_*(\Symp(\tMuc)) \to \pi_*(\Diff(\tM))
$
is generated by $e$.
\end{prop}

\subsection{Ample line bundles and Liouville vectors fields on some
algebraic open sets}\label{ss:Liouville}

In this subsection, we prove the following.
\begin{lemma}\label{le:Liouville}
Consider the manifold $(\tM, \tilde{\om}_{\mu, c})$ with $\mu$ rational. Let $i$ be any nonnegative integer and let $C_i \in {\mathcal C}_i$ be a configuration. Set $S_i = \tM - C_i$, the open subset of $\tM$ which is the complement of the union of the three spheres of $C_i$. Then there is on $S_i$ a Liouville vector field with only one zero and such that, for any compact set $K$ of $S_i$, there is a smooth hypersurface enclosing $K$ on which the field is transverse.
\end{lemma}

Recall that $C_i$ is defined in the following way. In the blown-up space $(\tM, \tilde{\om}_{\mu, c})$, let $D_{2k}, k \ge 0,$ denote the homology class $B - kF -E$ and let $C_{2k}$ be any configuration consisting of three embedded symplectic spheres, crossing transversally and positively, in classes $F-E, E, D_{2k}$. This is nonempty if and only if  $k + c < \mu$. Similarly,  let $D_{2k-1}, k \ge 1,$ denote the homology class $B - kF$, and let $C_{2k-1}$ be any configuration consisting of three embedded symplectic spheres, crossing transversely and positively, in classes $E, F-E, D_{2k-1}$ (this is non-empty
iff $k < \mu$). If $\Cc_i$  denote the corresponding configuration spaces, there is a natural ordering of these according to the decreasing area of the third sphere: $\Cc_0,\Cc_1,\Cc_2,\ldots$

\MS
\begin{proof}[Proof of Lemma~\ref{le:Liouville}]
In order to construct the Liouville field, one can in principle choose on $S_i$
coordinates in $\C^2$ which are holomorphic with respect to some natural projective structure on $S_i$ and then pull back the symplectic form. This would lead to the expression of the form as a potential. The problem is that the symplectic blow-up operation is not globally defined except on the projective space itself (by {\em globally defined} we mean that the blow-up $\tM$ of $M$ can be expressed as a K\"ahler submanifold of $M \times \CP^{n-1}$ which is a section over $M-\{pt\}$). We must therefore first approximate the symplectic form by a rational one; then embed the Hirzebruch surfaces (with rational forms) in a high-dimensional projective space using Kodaira's theorem and, finally, blow up that projective space and lift the Hirzebruch surface (i.e. take its proper transform).

Recall from Section~\ref{se:cohomology} that for any $\mu > 0$ and any integer $i\ge 0$ satisfying $\mu-\frac i2>0$, one may consider $\CP^1\times\CP^2$ endowed with the K\"ahler form $(\mu - \frac i2) \tau_1 + \tau_2$, where
$\tau_{\ell}$ is the Fubini-Study form on $\CP^{\ell}$ normalized so that the area of the linear $\CP^1$'s be equal to $1$. Let $W_i$ be the corresponding Hirzebruch surface, that is, the K\"ahler surface defined by
$$
W_i=\{([z_0,z_1],[w_0,w_1,w_2])\in\CP^1\times\CP^2 ~~|~~ z_0^i w_1=z_1^i w_0\}
$$

Now assume that the form $\om_\mu$ is rational on $S^2 \times S^2$, that is, that $\om_{\mu}=\mu \om(1) \oplus \om(1)$ with $\mu$ rational. Let $C_{2k} = (F-E), E, D_{2k}$ be a configuration as above in the space $\tMuc$, the blow-up of
$M_{\mu}$ at the ball of capacity $c < 1$ centered at the point $\bar p$. There is a K\"ahler structure on $M_{\mu}$ corresponding to $W_{2k}\subset
(\CP^1 \times \CP^2, (\mu - k)\tau_1 \oplus \tau_2)$. The curve $B-kF$ is represented by $s_0 = \{([z_0,z_1],[0,0,1])\}$, and we can choose the curve $F$ to be $\{([1,0],[w_0,0,w_2])\}$, so that the blow-up is centered at $\bar p = (B-kF) \cap F =([1,0],[0,0,1])$. Set $p/q = \mu - k$.  We can then realize
this structure projectively. This is what we first describe explicitly.

Consider the complete linear system of divisors of the ample line bundle
$(\pi_1^{*}H)^p \otimes (\pi_2^{*}H)^q$ over $\CP^1 \times \CP^2$ where $H$ denotes the hyperplane bundle over $\CP^1, \CP^2$ and $\pi_i$ are the projections on the two factors. By Kodaira's theorem, this defines an embedding
$$
\Phi_{p,q} : \CP^1 \times \CP^2 \to \CP^N
$$
given by
$$
\begin{array}{ccl}
([z_0,z_1],[w_0,w_1,w_2]) & \mapsto & [z_0^p w_0^q, z_0^p w_0^{q-1} w_1^1, \ldots, z_0^p
w_2^q,  \\
& & z_0^{p-1} z_1 w_0^q, z_0^{p-1} z_1 w_0^{q-1}w_1^1, \ldots, z_0^{p-1} z_1 w_2^q, \\
& &  \ldots \ldots \ldots \\
& & z_1^p w_0^q, z_1^p w_0^{q-1} w_1^1, \ldots, z_1^p w_1^q, \ldots, z_1^p w_2^q]
\end{array}
$$
Note that the pull-back $\Phi_{p,q}^{*} (\frac 1q \tau_N)$, restricted to
$W_{2k}$, is isotopic to the form $\om_{\mu}$. Blow up $\CP^N$ at $\Phi_{p,q}(\bar p)$ with weight $0 < c < 1$ so that $\widetilde{\CP^N}$ can be represented by the proper transform of $\CP^N$ inside
$(\CP^N\times\CP^{N-1},\frac{1}{q}\tau_N\oplus c\tau_{N-1})$. Taking now the proper transform of $\Phi_{p,q}(W_{2k})$ in
$(\CP^N\times\CP^{N-1},\frac{1}{q}\tau_N \oplus c \tau_{N-1})$ gives an algebraic model for $\tilde{M}_{\mu,c}$:
$$
\begin{array}{ccccc}
& & & & \CP^N \times \CP^{N-1} \\
& & & & \hspace{-.6cm} s \uparrow  \quad  \downarrow \\
W_{2k} & \subset & \CP^1 \times \CP^2 & \stackrel{\Phi_{p,q}}{- \hspace
{-.2cm}
\longrightarrow} &
\CP^N
\end{array}
$$
Let us parametrize the set $S_{2k}=\tMuc-((F-E)\cup E\cup(B-kF-E))$. In $W_{2k}$, a natural parametrization of $W_{2k}-(F \cup (B-kF))$ is:
$$
\begin{array}{cccc}
 & \C^2 & \to & \CP^1  \times \CP^2 \\
       & (a,b) & \mapsto & ([a,1],[a^{2k},1,b])
\end{array}
$$
Since
$$
\begin{array}{cl}
\Phi_{p,q}(\bar p) = \Phi_{p,q}([1,0],[0,0,1]) = & [0, 0, \ldots, 1,  \\
 & 0, 0, \ldots, 0, \\
 &  \ldots \ldots \ldots \\
 & 0, 0, \ldots,0]
\end{array}
$$
the parametrization of $S_{2k}$ is then
$$
\begin{array}{cccl}
\zeta: &  \C^2 &  \to & \CP^N  \times \CP^{N-1} \\
       & (a,b) & \mapsto & (\zeta_1(a,b), \zeta_2(a,b))
\end{array}
$$
where
$$
\begin{array}{ccl}
\zeta_1(a,b) & = & [a^{2kpq}, a^{2kp(q-1)}, \ldots , a^p b^q, \\
             &   & a^{2k(p-1)q}, a^{2k(p-1)(q-1)}, \ldots, a^{p-1}b^q, \\
             &   &  \ldots \ldots \\
             &   & a^{2kq}, a^{2k(q-1)}, \ldots,1, \ldots, b^q]
\end{array}
$$
and $\zeta_2(a,b)$ is defined like $\zeta_1(a,b)$ but with the last term of the first row deleted. The point here is that both functions have at least one
nonzero constant term. We can now pull back the form
$\frac{\tau_N}{q}\oplus c\tau_{N-1}$ by $\zeta$ to get $\zeta_1^{*}(\frac{\tau_N}{q}) + \zeta_2^{*}(c\tau_{N-1})$. Because
both maps are holomorphic and since
$\tau_{\ell} = \frac{i}{2\pi} \p\bar{\p}\log f_{\ell}$
with $f_{\ell} = \sum_0^{\ell} \|z_i\|^2$, we have
$$
\zeta^{*}(\frac{\tau_N}{q} \oplus c\tau_{N-1}) =
      \frac{i}{2\pi} \p \bar{\p} \log (g_1^{1/q}(a,b) g_2^c(a,b))
$$
where both $g_1$ and $g_2$ are of the form $1 +$ a polynomial of even powers in the variables $\|a\|, \|b\|$ without constant term. Taking the gradient of
$g=\log (g_1^{1/q}(a,b) g_2^c(a,b))$ with respect to the metric induced by the
pull-back form and the standard complex structure on $\C^2$ gives a Liouville vector field. One easily checks that $g$ it has only one critical point and that it exhausts $\C^2$. Hence, for any compact $K$, there is a level set of $g$ that contains $K$ in its interior and, by definition, the Liouville field is transverse to that hypersurface.

\medskip
This proves the Proposition for $G_{2k}$ when $\mu$ is rational.

\MS
The proof of the Proposition in the case of configurations in ${\mathcal C}_{2k-1}$, i.e. those of the form $E, (F-E),  B-kF$, is entirely similar, being understood that one works in this case with the blow-up of Hirzebruch surfaces of type $W_{2k-1}$.
\end{proof}

\noindent
{\bf Remark.}
We could not extend Lemma~\ref{le:Liouville} to irrational values of $\mu$. Here is why. Note first that the complement of a small tubular neighborhood of $C_i$ is clearly not convex, since for instance, the complement of a neighborhood of $E$, say, is symplectomorphic to the complement of a ball. Thus one cannot expect that the directions of a Liouville flow satisfying the conditions of the Lemma converge as $p \in \tM$ approaches a point in $C_i$.  On the other hand, the surgery that defines the symplectic blow-up in general is not compatible with the expression  of the potential for $c=0$; hence the embedding of the blow-up in projective spaces seems necessary.  Now, in order to pass from $\mu$ rational to $\mu$ irrational, we would consider $\mu_0$ irrational given, as well as $\phi \in G_i$ a symplectic diffeomorphism of $S_i = (S_i)_{\mu_0,c}$ with compact support. We would choose a small neighborhood $V_0$ of say the exceptional curve $A_2$ (in class $E$), disjoint from the support of $\phi$, then inflate slightly the form $\tilde{\om}_{\mu_0,c}$ in $V$ so that after renormalization the resulting form be in class $[\tilde{\om}_{\mu_1,c}]$ for some rational $\mu_1$ close to $\mu_0$. We would now apply our Lemma for the rational $\mu_1$, but we need that the resulting isotopy $\phi_t$ of $\phi$ have compact support outside $V_0$ (or more precisely, outside its image $V_1$ by the inflation + renormalization process) in order to take back that isotopy to the original setting $\mu = \mu_0$. Thus we need that ${\rm supp}(\phi)$ remain far enough from $V_1$. This can of course be achieved by choosing the neighborhood $V_0$ smaller, but then the rational approximation $\mu_1$ of $\mu_0$ becomes finer (with larger denominator) and the Liouville flow is wilder at infinity. What one needs is to control the behavior of the flow at infinity as the rational approximation becomes finer, which seems an intractable problem. Note finally that one natural way of establishing Lemma~\ref{le:Liouville} for irrational values of $\mu$ is, after all, to proceed as we did: prove the rational case, compactify, and establish Theorem~\ref{Thm:Stabilite} on stability and try to derive Lemma~\ref{le:Liouville} from the fact that the homotopy type of $\Symp(\tM)$ remains invariant. This will be taken up elsewhere.

\BS

\section{The rational homotopy groups of $\Emb_{\om}(c,\mu)$ }
\label{se:homotopy}

We are now in position to prove the main theorem of this paper, that
we recall for the convenience of the reader:
\begin{theorem*}[1.7]\label{thm:rational-homotopy2}
Fix any value $\mu\in (1,2]$. Then
\begin{enumerate}
\item[(1)]  the spaces $\Im \Emb_{\om}(c,\mu), c < \la= \mu-1$, are all
homotopically equivalent to $S^2\times S^2$ (they have therefore the homotopy type of a finite CW-complex);
\item[(2)] the rational homotopy groups of $\Im \Emb_{\om}(c,\mu), c \ge \la =
\mu -1$, vanish in all dimensions except in dimensions $2, 3$, and $4$, in which cases we have $\pi_2 =\Q^2$, $\pi_3 = \Q^3$, and $\pi_4 = \Q$ (they do not have the homotopy type of a finite CW-complex).
\end{enumerate}
\end{theorem*}
\begin{proof}  Let's begin by the first case. By Proposition~\ref{StructureCohomologie3Strates}, for $c < \la$, the groups $\Symp(\tMuc)$ and $\Symp(M_{\mu,c},B_c)$ have the same rational homotopy type as the stabilizer of a point
$\Symp_p(M_\mu)\subset\Symp(M_\mu)$. Consider the commutative diagram of homotopy fibrations:
$$
\begin{array}{ccccc}
\Symp(M_{\mu},B_c) & \hookrightarrow & \Symp(M_{\mu}) & \to &
\Im \Emb_{\om}(c,\mu) \\
\downarrow & & \downarrow \id & & \downarrow \\
\Symp_p(M_{\mu})  & \hookrightarrow & \Symp(M_{\mu}) & \to &
S^2 \times S^2
\end{array}
$$
where the first vertical arrow is the above homotopy equivalence and the
last one is induced by the restriction map on $\Emb_{\om}(c,\mu)$. Because the first two vertical arrows are homotopy equivalences, so is the last one.

\BS 
For the case $c \ge \la$, we first prove the following.
\begin{lemma} \label{le:j_1} The map $j_1: \pi_1(\Symp(\tMuc))
\otimes \Q \to\pi_1(\Symp(M_{\mu})) \otimes \Q$ in the homotopy sequence (3) of the introduction is surjective for all values of $\mu\in (1,2]$ and of $0<c<1$.
\end{lemma}
\begin{proof} This is a statement simpler than the formulas of Proposition~\ref{prop:invert}. It relies simply on the definition of the tori $\tT^2_i$ of section Section~\ref{se:cohomology}. Recall from Section~\ref{se:cohomology} that the  generator $a_3 = y_1$
of $\pi_1(\Symp(\tMuc))\otimes \Q$ can be described in the following way. Take the Hirzebruch integrable complex structure on $(S^2\times S^2,\om_{\mu})$ given by the projectivization of the complex bundle $L_2\oplus\C\to\CP^1$, where $L_2$ is the line bundle with divisor $2[pt]$. This is biholomorphic to the projection on the first factor of
$$
W_2 = \{([z_0,z_1], [w_0,w_1,w_2]) \in \CP^1 \times \CP^2 : z_0^2 w_1 - z_1^2 w_0 = 0 \}
\to \CP^1
$$
endowed with the restriction to $W_2$ of the K\"ahler form $\la\om_{\CP^1} +
\om_{\CP^2}$. The Hamiltonian $S^1$-action that rotates the fibers of
the line bundle lifts to a Hamiltonian $S^1$-action on the blow-up of $W_2$ at any of the fixed points of the action. This can be seen by cutting a triangle
containing the vertex representing the fixed point in the image of the moment map. By choosing the size of this triangle to correspond to the capacity $c$ and the fixed point to lie on the intersection of the $\infty$-section with the typical fiber, one gets the desired lift of $a\in \pi_1(\Symp(M_{\mu})) \otimes \Q$
to the element $a_3 = y_1 \in \pi_1(\Symp(\tMuc)) \otimes \Q$.
\end{proof}
\BS

Now the substitution of the elements of Theorems~\ref{thm:homotopy-M} and~\ref{thm:homotopy-tildeM} (2) in the exact rational homotopy sequence of
(3) of the introduction yields:
$$
\begin{array}{cccccccc}
& \ldots &  0 & \to &  \Q &
\stackrel{\rho_4}{\longrightarrow} &  \pi_4(\Emb_{\om}(c,\mu))\otimes
\Q & \\
& \stackrel{\delta_3}{\longrightarrow} & 0 &
\stackrel{j_3}{\longrightarrow} &  \Q^2 &
\stackrel{\rho_3}{\longrightarrow} & \pi_3(\Emb_{\om}(c,\mu)) \otimes
\Q & \qquad (5) \\
& \stackrel{\delta_2}{\longrightarrow} & \Q &
\stackrel{j_2}{\longrightarrow} & 0 &
\stackrel{\rho_2}{\longrightarrow} & \pi_2(\Emb_{\om}(c,\mu)) \otimes
\Q  & \\
& \stackrel{\delta_1}{\longrightarrow} & \Q^3
&\stackrel{j_1}{\longrightarrow} & \Q &
\stackrel{\rho_1}{\longrightarrow} & \pi_1(\Emb_{\om}(c,\mu))\otimes \Q &
\\
  & \longrightarrow & 0. &&&&&
\end{array}
$$

This gives $\pi_i(\Emb_{\om}(c,\mu)) \otimes \Q = 0$ for all $i\ge 5$,
$\pi_4(\Emb_{\om}(c,\mu))
\otimes \Q \simeq \Q$ and $\pi_3(\Emb_{\om}(c,\mu)) \otimes \Q
\simeq \Q^3$. Lemma~\ref{le:j_1} yields at once the ranks of the lower
homotopy groups. The fact that the space $\Im \Emb(c,\mu)$ does not have the homotopy type of a finite CW-complex for $c \ge \la$ is a consequence of the theory of minimal models and will be postponed to~\cite{Pinso}.
\end{proof}

\BS \BS


%


\begin{thebibliography}{99}
%
%
\bibitem{Ab}
Abreu, M., Topology of symplectomorphism groups of $S^2\times S^2$,
Invent. Math. 131 (1998), 1--23.
%
%
\bibitem{AM}
Abreu, M., McDuff, D., Topology of symplectomorphism groups of
rational ruled surfaces, J. Amer. Math. Soc. 13 (2000), 971--1009.
%
%
\bibitem{An}
Anjos, S., Homotopy type of Symplectomorphism Groups of $S^2\times S^2$,
Geometry and Topology 6 (2002), 195--218.
%
%
\bibitem{Au}
Audin, M., The topology of torus actions on symplectic manifolds, Progress in Math.  93,
Birkh\"auser, 1991.
%
%
\bibitem{CS}
Cartan, H., Serre, J.-P., Espaces fibr\'es et groupes d'homotopie II. Applications,
C. R. Acad. Sci. Paris 234 (1952), 393--395.
%
%
\bibitem{Gr}
Gromov, M., Pseudo holomorphic curves in symplectic manifolds,
Invent. Math. 82 (1985), 307--347.
%
%
\bibitem{HLS}
Hofer, H., Lizan, V., Sikorav, J.-C., On genericity for holomorphic curves
in four-dimensional almost complex manifolds ,
J. Geom. Anal. 7 (1997), 149--159.
%
%
\bibitem{L}
Lalonde, F., Isotopy of symplectic balls, Gromov's radius,
and the structure of irrational ruled symplectic $4$\/-manifolds,
Mathematische Annalen 300 (1994), 273-296.
%
%
\bibitem{LMD}
Lalonde, F., McDuff, D., $J$-curves and the classification of rational
and ruled symplectic $4$-manifolds,
in: Contact and symplectic geometry, Publ. Newton Inst., 8, Cambridge
Univ. Press, Cambridge, 1996, pp 3--42.
%
%
\bibitem{LL}
Li T. J., Liu A., General wall crossing formula,
Math. Res. Lett. 2 (1995), 797-810.
%
%
\bibitem{MD:BlowUps}
McDuff, D., Blow ups and symplectic embeddings in dimension $4$,
Topology 30 (1991), 409--421.
%
%
\bibitem{MD:Isotopie}
McDuff, D., From Symplectic Deformation to Isotopy,
in {\em Topics in symplectic $4$-manifolds, Irvine, CA, 1996},
Internat. Press, Cambridge, 1998, pp 85--99.
%
%
\bibitem{MD:Stratification}
McDuff, D., almost complex structures on $S^2\times S^2$,
Duke Math. J. 101 (2000), 135-177.
%
%
\bibitem{MD:Inflation}
McDuff, D., Symplectomorphism groups and almost complex structures,
to appear in {\it Enseignement Math}.
%
%
\bibitem{MP}
McDuff, D., Polterovich, L., Symplectic packings and algebraic geometry,
Invent. Math. 115 (1994), 405--429.

%
%
\bibitem{Me}
Meigniez, G., Submersions, Fibrations and Bundles,
Trans. Amer. Math. Soc. 354 No. 9 (2002), 3771--3787.
%
%
\bibitem{Pinso}
Pinsonnault, M., In preparation
%
%
\bibitem{Ta}
Taubes, C. H., Counting pseudo-holomorphic submanifolds in dimension 4,
J. Differential Geom. 44 (1996), 818--893.
%
%
\end{thebibliography}
\end{document}